\documentclass[10pt]{amsart}
\usepackage{amscd,amsmath,amssymb}
\usepackage{hyperref}
\usepackage[all]{xy}
\usepackage{graphicx}
\usepackage{psfrag}

\newcommand{\al}{\alpha}

\newcommand{\Th}{\Theta}
\newcommand{\Ad}{\operatorname{Ad}}

\newcommand{\ov}{\overline}

\newcommand{\TT}{\mathbb{T}}

\newcommand{\RR}{\mathbb{R}}

\newcommand{\restr}[1]{\vrule height3ex width.4pt depth1.4ex\lower1.4ex\hbox{\scriptsize $\,#1$}}
\newcommand{\rrestr}[1]{\vrule height2ex width.4pt depth0.9ex\lower0.9ex\hbox{\scriptsize $\,#1$}}
\newcommand{\g}{\mathfrak{g}}

\numberwithin{equation}{section}

\newtheorem{te}{Theorem}[section]

\newtheorem{pr}{Proposition}[section]

\newtheorem{lm}{Lemma}[section]

\theoremstyle{definition}
\newtheorem{de}{Definition}[section]
\newtheorem{re}{Remark}[section]

\begin{document}

\title{Singular Cosphere Bundle Reduction }
\thanks{The first and the second authors thank the Swiss National Science Foundation
for partial support.}
\author{Oana M. Dr\u agulete }
\address{O. M. Dr\u agulete: Section de math\'ematiques, EPFL, CH-1015 Lausanne, Switzerland
\newline\indent
and
\newline\indent
D{e}partement of Math{e}matics, University ``Politehnica'' of
Bucharest, Romania} \email{oana.dragulete@epfl.ch}

\author{Tudor S. Ratiu}
\address{T. S. Ratiu: Section de math\'ematiques, EPFL, CH-1015 Lausanne, Switzerland}
\email{tudor.ratiu@epfl.ch}

\author{Miguel Rodr\'{i}guez-Olmos}
\address{M. Rodr\'{i}guez-Olmos: Section de math\'ematiques, EPFL, CH-1015 Lausanne, Switzerland} \email{miguel.rodriguez@epfl.ch}


\begin{abstract}
This paper studies singular contact reduction for cosphere bundles at
the zero value of the momentum map. A stratification of the singular
quotient, finer than the contact one and better adapted to the bundle
structure of the problem, is obtained. The strata of this new
stratification are a collection of cosphere bundles and coisotropic
or Legendrian submanifolds of their corresponding contact components.
\end{abstract}

\maketitle

\section{Introduction}
The main goal of this paper is to carry out the  singular
reduction of cosphere bundles at the zero value of the contact
momentum map. This presents interest because cosphere bundles
carry considerably more structure than a general contact manifold
and have the same privileged position in contact geometry that
cotangent bundles have in symplectic geometry.

Contact reduction appears for the first time in the work of
Guillemin and Sternberg \cite{gs} in the context of reducing
symplectic cones. Albert \cite{al} and, several years later,
Geiges \cite{ge} and Loose \cite{loose}  independently defined and
studied contact reduction at the zero value of the contact
momentum map for free proper contact actions of Lie  groups.
Reduction at a general value of the momentum map was studied by
both  Albert \cite{al} and Willett \cite{wil} who proposed two
different versions of dealing with it. It turns out that Willett's
method is the one that naturally parallels the symplectic
reduction theory, even in the singular case as shown by Lerman and
Willett \cite{lw}. They prove that the resulting  contact
quotient, at any value of the momentum map, depends  only on the
contact structure, that it is independent of any contact form that
defines the contact foliation, and that it is a stratified space,
more precisely, a cone space.

The case of cosphere bundle reduction for proper free lifted Lie
group actions was studied in \cite{dor} with a view of comparing
the  theory to that for cotangent bundle reduction. It turns out
that in regular contact reduction of cosphere bundles there are no
analogues of magnetic terms. In parallel, in \cite{prs} the
authors have developed the theory of singular cotangent bundle
reduction at the zero value of the momentum map and have found a
finer stratification than that given by the general theory due to
the additional structure  of the cotangent bundle and the fact
that the Lie group action is a cotangent lifted action. A similar
phenomenon occurs in contact reduction of cosphere  bundles.
Applying the general theory of singular contact reduction due to
Lerman and Willett \cite{lw} yields contact stratified  spaces
that, however, lose all information of the internal structure of
the cosphere bundle. Based on the cotangent bundle reduction
theorems, both in the regular and singular case, as well as
regular cosphere bundle  reduction, one expects additional
bundle-like structure for the contact strata. The cosphere  bundle
projection to the base manifold descends  to a continuous
surjective map from the reduced space at zero to the orbit
quotient of configuration space, but it fails to be a  morphism of
stratified spaces if we endow the reduced space with its contact
stratification and the base space with the customary orbit type
stratification defined  by the Lie group action. The present paper
introduces a new stratification of the contact quotient at zero,
called in what  follows the \textit{C-L stratification} (standing
for the coisotropic or Legendrian nature of its pieces) which
solves the above mentioned two problems. Its main features are the
following. First, it is compatible with the contact stratification
of the quotient and the orbit type stratification of the
configuration orbit space. It is also finer than the contact
stratification. Second, the natural projection of the C-L
stratified quotient space to its base space, stratified by orbit
types, is a morphism of stratified spaces. Third, each C-L stratum
is a bundle over an orbit type stratum of the base and each
contact stratum can be seen as a union of C-L pieces, one of them
being open and dense in its corresponding contact stratum and
contactomorphic to a cosphere bundle. The other strata are
coisotropic or Legendrian submanifolds  in the contact components
that  contain them.

The paper is structured as follows. Section 2 presents the
definitions, conventions, and results on stratified spaces and
contact reduction (regular and  singular) that are used throughout
the paper. Section 3 quickly reviews the relevant results on
regular contact cosphere reduction. Section 4 presents the
stratification of  the zero level set of the momentum map and
begins  the  work on the stratification of the quotient by
studying the case of one single orbit type (Theorem \ref{uf}). The
contact stratification and contact geometry of the reduced space
are studied in Section 5, having as main results Theorems
\ref{theoremsecondary} and \ref{theoremsecondaryb}. The new C-L
stratification is also introduced here and its properties are
investigated. Theorem \ref{coisotropicstratification} presents a
complete description of its frontier conditions. Section 6 studies
the singular cosphere bundle reduction for almost semifree
actions, that is, actions that are in bijective correspondence
with free lifted actions on the cosphere bundle. The
stratification is computed explicitly and the particular case of
the circle acting on the cosphere bundle of the plane is carried
out in detail. Section 7 studies the example of the diagonal
action of the two-torus on two copies of the plane, lifted to the
cosphere bundle. This example is rich enough to illustrate the
relationships between the various stratifications and the strata
are computed explicitly.

\section{Preliminaries}
In this section we will survey the main results of several topics
that will be needed in the subsequent development of the paper. We
will assume that all topological spaces are paracompact. In
addition, manifolds will be real, smooth and finite-dimensional.
By group we will mean a finite-dimensional Lie group. Every action
of a group $G$ on a manifold $M$ is supposed to be smooth and the
usual notation $g\cdot m$ for $g\in G$ and $m\in M$ will be
employed. The natural pairing between a vector space and its dual
will be denoted by $\langle \cdot, \cdot \rangle$. By submanifold,
we will always mean an embedded submanifold.
\subsection{Stratified spaces and proper group actions}
The natural framework for singular reduction is the category of
stratified spaces. We briefly recall here the basic concepts (see
\cite{pf}). Let $X$ be a topological space and
$\mathcal{Z}_X=\{S_i:i\in I\}$ a locally finite partition of $X$
into locally closed  disjoint subspaces $S_i\subset X$, where $I$
is some index set. We  say that $(X,\mathcal{Z}_X)$ is a
\textit{decomposed space} if every  $S_i$ is a manifold whose
topology coincides with the induced one from $X$ and if the
\textit{frontier condition} holds:
$S_i\cap\overline{S_j}\neq\varnothing$ implies
$S_i\subset\overline{S_j}$, whence  $S_i\subset
\partial S_j$, where
$\partial S_j:=\overline{S_j}\backslash S_j$. In this case, the
elements of $\mathcal{Z}_X$ are called \textit{pieces\/} of the
decomposition.

In a topological space $X $, two subsets $A $ and $B $ are said to
be \textit{equivalent at\/} $x $ if there exists an open
neighborhood $U $ of $x $ such that $A \cap U = B \cap U $. These
equivalence classes are called \textit{set germs at} $x $. Let
$\mathcal{S}$ be the map that associates to each point $x\in X$
the set germ $\mathcal{S}_x=\left[O\right]_x$ of a locally closed
subset $O$ of $X$. We say that $(X,\mathcal{S})$ is a
\textit{stratified space\/} if, for every point $x\in A$ there
exists a neighborhood $U$ of $x$ endowed with a decomposition
$\mathcal{Z}_U$ such that for every $y\in U$,
$\mathcal{S}_y=\left[Z(y)\right]_y$, where $Z(y)\in\mathcal{Z}_U$
denotes the piece containing $y$. In this case we say that the
decomposition $\mathcal{Z}_U$ \textit{locally induces\/}
$\mathcal{S}$.

Given two stratified spaces $(X,\mathcal{S})$ and
$(Y,\mathcal{T})$ and a continuous map $f:X\rightarrow Y$, we say
that $f$ is a \textit{morphism of stratified spaces\/} (or shorter, a
\textit{morphism}) if for every $x\in X$ there exist neighborhoods $V$
of
$f(x)$ and $U\subset f^{-1}(V)$ of $x$ such that
\begin{itemize}
\item[(i)] there exist decompositions $\mathcal{Z}_U$ and
$\mathcal{Z}_V$ locally inducing the stratifications $\mathcal{S}$
and $\mathcal{T}$ respectively, with the property that for every
$y\in U$ contained in a piece $S\in\mathcal{Z}_U$ there is an open
neighborhood $y\in W\subset U$ such that $f\rrestr{W} (S\cap W)$
is contained in the unique piece $R\in\mathcal{Z}_V$ that contains
$f(y)$, and \item[(ii)] $f\rrestr{S\cap W}:S\cap W\rightarrow R$
is smooth.
\end{itemize}
In addition, we will say that $f$ is a \textit{stratified
immersion\/} (resp. \textit{submersion, diffeomorphism}, etc...)
if so are all the maps $f\rrestr{S\cap W}$ for every point $x\in
X$. Given two different stratifications $\mathcal{S}$ and
$\mathcal{S}'$ on the same topological space $X$, we say that
$\mathcal{S}$ is \textit{finer\/} than $\mathcal{S'}$ if the
identity map $1_X$, viewed as a map between stratified spaces
$(X,\mathcal{S})\rightarrow (X,\mathcal{S}')$, is a morphism.

Smooth manifolds are trivially stratified spaces and smooth maps
between manifolds are their morphisms. Note that a decomposed
space $(X,\mathcal{Z}_X)$ induces naturally a stratification
$(X,\mathcal{S})$ by just taking $S_x$ to be the set germ of the
piece containing $x$, for every $x\in X$. In this case, we call
the pieces $S_i\in\mathcal{Z}_X$ the \textit{strata} of
$(X,\mathcal{S})$ and say that they satisfy the frontier
conditions defined by the underlying decomposition. In this paper
the stratifications that will appear will be of this form and
thus, for the sake of simplicity, when this is the case we will
work most of the time with the decompositions inducing these
stratifications.
\medskip

Let $\phi:G\times M\rightarrow M$ be a smooth action of the
Lie group $G$ on the  manifold $M$. Since $M $ is paracompact it
admits a Riemannian metric so if it is connected, $M $ is second
countable. The action is called \textit{proper\/} if $\phi \times
\operatorname{id}_M$ is a proper map. In this paper we only work
with proper actions. For instance, every action of a compact group
is automatically proper. The main properties of a proper action of
$G$ on $M$ are:\begin{itemize} \item[(i)] For each $m\in M$, its
stabilizer (or isotropy group) $G_m$ is compact. \item[(ii)] The
manifold structure of the orbit $G\cdot m$ is the one that makes
the natural bijection $G/G_m \rightarrow G \cdot m $ a
diffeomorphism. The inclusion $G \cdot m \hookrightarrow M$ is an
injective immersion. In addition, the orbit is a closed subset of
$M $. If $M $ is connected, then the orbit is an embedded
submanifold of $M $. \item[(iii)] The quotient space equipped with
the quotient topology is  paracompact and the orbit map
$\pi:M\rightarrow M/G$ is open and closed. \item[(iv)] $M$ admits
a $G$-invariant Riemannian metric. \item[(v)] If all the
stabilizer groups are conjugate to a given subgroup $H\subset G$,
then $M/G$ is a smooth manifold, the orbit map
$\pi:M\rightarrow M/G$ is a smooth locally trivial fiber bundle
whose fibers are diffeomorphic to $G/H$, and the structure group
of this locally trivial fiber bundle is $N(H)/H $, where $N(H)$
is the normalizer of $H$ in $G$.
\end{itemize}

We now quote Palais' Tube Theorem \cite{Palais} in a form adapted
to our needs, which is of great importance in the local study of
proper actions. Let $m\in M$. Choose an invariant Riemannian
metric on $M$ and use it to decompose $T_mM=\mathfrak{g}\cdot
m\oplus S_m $, where $\g\cdot m=\{\xi_M(m)\, :\, \xi\in\g\}$. This
splitting is
$G_m$-invariant for the linear action of $G_m$ on $T_mM$. The
twisted action of $G_m$ on
$G\times S_m$ is defined by
\begin{equation}\label{twist}
h\cdot (g,s)=(gh^{-1},h\cdot s)
\end{equation}
for $h\in G_m, g\in G$ and $s\in S_m$. Since $G_m$ acts freely on
the right on $G$, the twisted action is free. In addition, $G_m$
is compact by property (i) of proper group actions, so the
quotient space, denoted by $G\times_{G_m}S_m$, is a manifold. The
Tube Theorem implies the existence of a $G_m$-invariant open ball
$U$ around the origin in $S_m$ such that the map
$\psi:G\times_{G_m}S_m\rightarrow M$ defined by
\begin{equation}\label{tube}
\psi ([g,s])=g\cdot \mathrm{exp}_m (s)
\end{equation}
maps $G\times_{G_m}U$ diffeomorphically and equivariantly onto a
$G$-invariant neighborhood $U'$ of $G\cdot m$ in $M$. Here,
$\mathrm{exp}_m$ is the exponential map at $m$ associated to the
chosen Riemannian metric. The map $\psi$ is called a
\textit{tube\/} for the action and $S_m$ is called a
\textit{linear slice\/}, or simply a \textit{slice\/} of the
action at $m $.

Let $I_M$ be the isotropy lattice of $M$, i.e. the set of
conjugacy classes of subgroups of $G$ which appear as stabilizers
for the action of $G$ on $M$. Such classes, called \textit{orbit
types}, are denoted by $(H) $ . For each element $(H)\in I_M$ the
$(H)$-orbit type manifold is defined by
\begin{equation}\label{orbit type}
M_{(H)}=\{m\in M \mid (G_m)=(H)\}.
\end{equation}
In the same way, for any subset $A$ of $M$ one defines the orbit
type sets of $A$ by $A_{(H)}=A\cap M_{(H)}$ and the isotropy
lattice of $A$ by restriction. For a proper $G$-action on a manifold $M$ such
that $M/G$ is connected, there is always a subgroup $H_0\subset G$
such that $M_{(H_0)}$ is open and dense in $M$ and $H_0$ is
conjugate to a proper subgroup of any other stabilizer. This
orbit type $(H_0)$ is called the \textit{principal orbit type\/}
of $I_M$.

Obviously, the collection of orbit type manifolds forms a partition
of $M$. For simplicity, we will make from now on the following
important assumption: for every $(H)\in I_M$, all the connected
components of $M_{(H)}$ have the same dimension and $M$ is second
countable. Hence we have:
\begin{itemize}\item[(i)] For every $(H)\in I_M$, $M_{(H)}$ is
a $G$-invariant submanifold of $M$, and \item[(ii)] $M$ and $M/G$
are stratified spaces with strata $M_{(H)}$ and
$M^{(H)}:=M_{(H)}/G$ respectively. Their frontier conditions are:
$$M^{(H)}\subset\partial M^{(L)}\Longleftrightarrow (L)\prec
(H),$$ and correspondingly for $M$, where $(L)\prec (H)$ means
that $L$ is conjugate to a proper subgroup of $H$. Since $\prec$
defines a partial ordering in $I_M$ we say that the frontier
conditions of the stratification of $M/G$ are induced by the
isotropy lattice $I_M$.
\end{itemize}

\begin{re}
If one allows the connected components of the orbit type manifolds
to have different dimensions, then one needs to work in the larger
category of $\Sigma$-manifolds and $\Sigma$-decompositions. A
$\Sigma$-manifold is a countable topological sum of connected
smooth manifolds having possibly different dimensions (see \cite{pf} for more details). However, our
results on the stratified nature of the studied quotient spaces
remain valid.
\end{re}

\subsection{Reduction of Contact Manifolds.}
 Recall that a \textit{contact structure} on a smooth
$(2n+1)$-dimensional manifold $\mathcal{C}$  is a codimension one
smooth distribution $\mathcal{H} \subset T\mathcal{C}$ maximally
non-integrable in the sense that it is locally given by the kernel
of a one-form $\eta$ with $\eta \wedge (d\eta)^n\ne 0$. Such an
$\eta$ is called a (local) \textit{contact form}. Any two
proportional contact forms define the same contact structure. A
contact structure which is the kernel of a global contact form is
called {\em exact}. In the case of exact contact manifolds,
$d\eta$ has rank $n$ implying the existence of the Reeb vector
field $R$ uniquely defined by
$$
i_Rd\,\eta=0 \quad\text{and}\quad \eta(R)=1.
$$
In the following we will consider only exact orientable contact
manifolds.

When studying the geometry of the singular reduced spaces of
cosphere bundles one needs the notions of coisotropic and
isotropic submanifolds in the contact context. Any integral
submanifold $N$ of $\mathcal{H}$ has the property that
its tangent space at every
point is an isotropic subspace of the symplectic vector space
$(\ker\,\eta_x,d\,\eta_x)$ and that's why, sometimes, they are
also called \textit{isotropic} submanifolds. In particular,
$\operatorname{dim}N\leq n$; if $\operatorname{dim}N=n$, then $N$
is called a \textit{Legendrian} submanifold. A submanifold $N$ of
the  contact manifold $(\mathcal{C},\eta, R)$ is
\textit{coisotropic\/} if for any $x\in N$ the subspace $T_x N
\cap \ker \eta_x $ is coisotropic in the symplectic vector space
$(\ker\,\eta_x,d\,\eta_x)$.

 A group $G$ is said to act by
\textit{contactomorphisms} on a contact manifold if it preserves
the contact structure $\mathcal{H}$. For an exact contact manifold
$(\mathcal{C},\eta)$, this means that $g^*\eta=f_g \eta$ for a
smooth, real-valued, nowhere zero function $f_g$. $G$ acts by
\textit{strong contactomorphisms\/} on $\mathcal{C}$, if
$g^*\eta=\eta$, i.e. $G$ preserves the contact form, not only the
contact structure. A $G$--action by  strong contactomorphisms on
$(\mathcal{C}, \eta)$ admits an equivariant momentum map
$J:N\rightarrow \mathfrak{g}^*$ given by evaluating the contact
form on the infinitesimal generators of the action: $\langle J
(x), \xi\rangle := \eta(\xi_{\mathcal{C}})(x)$. Note the main
difference with respect to the symplectic case: any action by
strong contactomorphisms is automatically Hamiltonian. Note also
that orbits which lie in the zero level set of the contact
momentum map are examples of isotropic submanifolds.
 For more details on contact manifolds and their associated momentum
 maps see \cite{bl}, \cite{ge1}, and \cite{wil}.

Reduction theory for co-oriented contact manifolds in the singular
context was introduced by Willett
 in \cite{wil}. We now review briefly this construction at zero momentum,
 since it will be used in our next refinement to the cosphere bundle case.
  Let $G$ be a  group that acts
 by strong contactomorphisms on an exact contact manifold
$(\mathcal{C},\eta)$. By the definition of the momentum map, its
zero level set is a $G$-space. The contact quotient (reduced
space) of $\mathcal{C}$ at zero momentum is defined as
$$
\mathcal{C}_0 : = J^{-1} (0)/G.
$$
Note that, as in the symplectic case, this quotient is in general
a singular space.

\begin{te}
\label{WL1} Let $(\mathcal C,\eta)$ be an exact contact manifold and
$G$ a Lie group acting properly on $\mathcal C$ by strong
contactomorphisms.
 Then for every stabilizer subgroup $H$ of $G$ the set
$$
\mathcal{C}_0^{(H)}:=(J^{-1}(0))_{(H)}/G=(\mathcal{C}_{(H)}\cap
J^{-1}(0))/G
$$
is a smooth manifold  and the partition of the contact quotient
$$
\mathcal C_0:=\left(J^{-1}(0)\right)/G
$$
into these manifolds is a stratification with frontier condition
induced by the partial order of $I_{J^{-1}(0)}$. Moreover, there
is a reduced exact contact structure on $\mathcal{C}_0^{(H)}$
 generated by the one-form $\eta_0^{(H)}$ characterized by
$$
\pi_G^{(H)}\eta_0^{(H)}=\tilde{i}_{(H)}\eta,
$$
where $\pi_G^{(H)}:(J^{-1}(0))_{(H)}\rightarrow
\mathcal{C}_0^{(H)}$ is the projection on the orbit space and
$\tilde{i}_{(H)}:(J^{-1}(0))_{(H)}\hookrightarrow \mathcal{C}$ is
the inclusion.
\end{te}

In what follows this stratification will be referred to as the
{\textit{contact stratification} of $\mathcal C_0$.

\section{Regular cosphere bundle reduction}
Cosphere bundles are the odd dimensional analogs of cotangent
bundles in contact geometry. In the following, we will briefly
recall their construction and their equivariant regular contact
reduction, referring to \cite{dor} and \cite{rs} for more details.

Let $Q$ be a $n$-dimensional manifold and $\theta$ the Liouville
one-form on $T^*Q$, defined by $\theta (X_{p_x})=\langle
p_x,T_{p_x}\tau X_{p_x}\rangle$, where $p_x\in T^*_xQ$, $X\in
T_{p_x}(T^*Q)$, and $\tau :T^*Q\rightarrow Q$ is the canonical
projection. Let $\Phi: G \times Q \rightarrow Q$ be an action of
 $G$ on $Q$. Denote by
\[
\Phi_*:G\times T^*Q\rightarrow T^*Q
\]
 its natural (left) lift to the cotangent
bundle. Consider  the action of the multiplicative group
 $\RR_+ $ by
dilations on the fibers of  $T^*Q\setminus\{0_{T^*Q}\}$.
\begin{de}
The {\em  cosphere bundle} $S^\ast Q$ of $Q$ is the quotient manifold
$(T^*Q\setminus\{0_{T^*Q}\})/\RR_+$.
\end{de}
Let $\pi_+:T^*Q\setminus\{0_{T^*Q}\}\rightarrow S^*Q$ and $\kappa:
[\alpha_q] \in S^\ast Q \mapsto q\in Q$ be the canonical
projections. Denote by $[\al_q]$ the elements of the cosphere
bundle. Of course, $(\pi_+, \RR_+, T^*Q \setminus\{0_{T^*Q}\},
S^*Q)$ is a $\RR_+$-principal bundle. Also, we will use the
$\pi_+$ notation for any $\RR_+$ projection. The exact contact
structure of $S^\ast Q$ is given by the kernel of any one form
$\theta_\sigma$ satisfying $\theta_\sigma=\sigma^*\theta$ for
$\sigma:S^*Q\rightarrow T^*Q\setminus\{0_{T^*Q}\}$ a global
section. Such $\sigma$ always exists and, even more, the set of
global sections of this principal bundle is in bijective
correspondence with the set of $C^\infty$ functions
$f:T^*Q\setminus\{0_{T^*Q}\}\rightarrow \RR_+$ satisfying
\[
f_\sigma(r \al_q)=\displaystyle\frac 1r  f_\sigma(\al_q),
\quad r\in \RR_+, \; \al_q\in T^*Q\setminus\{0_{T^*Q}\}.
\]
(See \cite{dor} for details).

\begin{re}
1. Let $\mathcal{C}(S^*Q)=S^*Q\times \RR_+$ be the symplectic cone
over $S^*Q$, endowed with the symplectic form $d(t\theta_\sigma)$.
Then one can easily see that
$T_\sigma:\mathcal{C}(S^*Q)\rightarrow T^*Q$ given by
$T_\sigma([\al_q], t)=tf_\sigma(\al_q) \al_q$ is a well
defined symplectic diffeomorphism, that is, a symplectomorphism.\\
2. If $Q$ is zero-dimensional, we set, by convention, $S^*Q=\varnothing$.
\end{re}
The action $\Phi$ lifts to the cosphere bundle yielding a proper
action
\[
\widehat \Phi_*:G\times S^*Q\rightarrow S^*Q \text{,}\quad
\widehat\Phi_*(g,[\al_q])=[\Phi_*(g,\al_q)]
\]
by contactomorphisms with all scale factors positive. In \cite{le}
it has been proved that for any proper action which preserves an
exact contact structure, there exists a $G$-invariant contact
form. As every contact form on the cosphere bundle is obtained
$via$ a global section as above,  we shall chose once and for all
a section $\sigma$ for which
$(\widehat\Phi_{*g})^*\theta_\sigma=\theta_\sigma$. Relative to
this contact form the induced action on the cosphere bundle is by
strong contactomorphisms. The associated momentum map, which
depends on the section $\sigma$, will be denoted $J$ for
simplicity, since in what follows no other contact form different
from $\theta_\sigma$ will be used. As above, the exact contact
structure of $S^*(Q/G)$ can be described as the kernel of a global
contact form of type $\Theta_\Sigma$, where
\[
\Sigma: S^*(Q/G)\rightarrow T^*(Q/G)\setminus \{0_{T^*(Q/G)}\}
\]
is a global section, and $\Theta$ is the Liouville one-form of
$T^*(Q/G)$.

Regular reduction of cosphere bundles was done in \cite{dor}. Its
main result at zero momentum is
\begin{te}
\label{olt1} Let $G$ be a finite dimensional Lie group acting
freely and properly on a differentiable manifold $Q$. Then
$(S^*Q)_0$, the reduced space at the regular value zero of the
cosphere bundle of $Q$, is contact diffeomorphic to the cosphere
bundle $S^{*}(Q/G)$.
\end{te}

In the remainder of this paper, we will generalize this result to
non-free actions, within the framework of stratified spaces,
relating our results to the contact stratification defined in
Theorem \ref{WL1}.

\section{The decomposition of $J^{-1}(0)$}

The geometric study of the contact reduced space $(S^*Q)_0$ passes
through the analysis of the level set $J^{-1}(0)$ and, in
particular, of its isotropy lattice $I_{J^{-1}(0)}$. We shall use
the fact that both the cosphere bundle $S^*Q$ and the lifted
action of $G$ on it are completely determined by the differential
structure of $Q$ and its supported $G$-action. This will allow us
to obtain our first main result, Proposition \ref{propofirst},
which describes this isotropy lattice, and hence the topology of
the contact stratification of $(S^*Q)_0$, in terms of the isotropy
lattice of $Q$ without those elements corresponding to zero
dimensional orbit types. Also, as a preliminary result, and a
``building block'' for the general construction, we state an
intermediary cosphere reduction result, Theorem \ref{uf}, which
applies to base manifolds $Q$ on which the group action is not
free but exhibits a single orbit type, that is, $I_Q$ consists of
only one element.

\begin{lm}
\label{latice}
 The isotropy lattice of the cosphere bundle coincides with
the isotropy lattice of the cotangent bundle without the zero
section
$$
I_{S^*Q}=I_{T^{*}Q\setminus \{0_{T ^\ast Q}\}}.
$$

\begin{proof}
It is enough to show that $G_{\al_q}=G_{[\al_{q}]}$ for any
$\al_q\in T^{*} Q\setminus\{0_{T ^\ast Q}\}$. Thus let $g\in
G_{[\al_{q}]}$. This implies that $g[\al_{q}]=[g\al_{q}]=[\al_{q}]
\iff g\al_{q}=r\al_{q}$ for $r>0$. Since the action of $G$ on $Q$
is proper, there is a $G$-invariant Riemannian metric on $Q$ and
hence $\Vert g\al_{q} \Vert=\Vert\al_q \Vert=r\Vert\al_q \Vert$.
It follows that $r=1$ and $G_{[\al_ {q}]}\subset G_{\al_q}$. The
other inclusion being obvious, the proof is now complete.
\end{proof}
\end{lm}

\begin{re}
We will write $J_{ct}:T^*Q\rightarrow\mathfrak{g}^*$ for the
canonical momentum map for the cotangent-lifted action of $G$ on
$T^*Q$ endowed with the canonical symplectic form. As
$J^{-1}(0)=\pi_+ (J_{ct}^{-1}(0)\setminus \{0_{T ^\ast Q}\})$ note
that

$$(J^{-1}(0))_{(L)}=\pi_+ \left((J_{ct}^{-1}(0))_{(L)}
\setminus [ (J_{ct}^{-1}(0))_{(L)} \cap
\{0_{T ^\ast Q}\}]\right)$$
\end{re}
since $(J_{ct}^{-1}(0)\setminus \{0_{T^\ast Q}\})_{(L)}
= (J_{ct}^{-1}(0))_{(L)}
\setminus [ (J_{ct}^{-1}(0))_{(L)} \cap \{0_{T ^\ast Q}\}]$.
\medskip

The following theorem is an immediate consequence of Theorems
$\ref{olt1}$ and $\ref{WL1}$.

 \begin{te}
\label{uf} Let $G$ be a finite dimensional Lie group acting
properly on the differentiable manifold $Q$ such that all the
points in $Q$ have stabilizers conjugate to some $K$ {\rm (}that
is, $Q=Q_{(K)}${\rm )}. Then $J^{-1}(0)$ is a submanifold of
$(S^*Q)_{(K)}$ and $(S^*Q)_0$, the reduced space at zero, is
contact-diffeomorphic to $S^{*}(Q/G)$ .
\end{te}

In the following proposition we give the decomposition of
$J^{-1}(0)$ and show how the topology of the contact quotient at
zero is completely determined by the isotropy lattice of $Q$. For
that, we will use the following partition of $T^*Q$. We fix once
and for all a $G$-invariant Riemannian metric on $Q$. Then, for
any $(H)\in I_Q$, the restriction of $TQ$ to the embedded
submanifold $Q_{(H)}$ can be decomposed as the Whitney sum
$T_{Q_{(H)}}Q=TQ_{(H)}\oplus NQ_{(H)}$, where, for every $q\in
Q_{(H)}$, $N_qQ_{(H)}=T_qQ_{(H))}^\perp$. Note that each of the
elements of the Whitney sum are $G$-invariant vector bundles over
$Q_{(H)}$. Dualizing this splitting over each orbit type
submanifold in $Q$, we obtain the following $G$-invariant
partition of $T^*Q$:
$$T^*Q=\coprod_{(H)\in I_Q}T^*Q_{(H)}\oplus N^*Q_{(H)}.$$
Now, the restriction of this partition to $T^*Q\backslash\{
0_{T^*Q}\}$ and afterwards its quotient by the action of
$\mathbb{R}_+$, induces a $G$-invariant partition of $S^*Q$.

Let $I^*_Q$ denote the isotropy lattice of $Q$ without those
elements corresponding to orbit type manifolds $(H)$ for which the
orbits of the restricted $G$-action have the same dimension as
$Q_{(H)}$. At this moment, we will need some results on
cotangent-lifted actions, which were proved in \cite{prs}.
\begin{lm}\label{lemmacotlift} If $G$ acts on $Q$ and on $T^*Q$ by cotangent lifts
with momentum map $J_{ct}:T^*Q\rightarrow \mathfrak{g}^*$. Let
$(L),\,(H)\in I_Q$ be arbitrary.
\begin{enumerate}
\item[(i)]
 $(N^*Q_{(H)})_{(H)}$ is the zero section
of $N^*Q_{(H)}$.
  \item[(ii)] Let $J_{ct(H)}$ denote the canonical
momentum map on $T^*Q_{(H)}$ associated to the lift of the action
on $Q_{(H)}$ obtained by restriction from $Q$. Then
\begin{equation}
\label{descomp}
(J_{ct}^{-1}(0))_{(L)}=J^{-1}_{ct(L)}(0)\coprod_{(H)\succ(L)} \left(
J^{-1}_{ct(H)}(0)\times (N^{*}Q_{(H)})_{(L)} \right)
\end{equation}
 \item[(iii)] If $(L)\neq (H)$, then
$(N^*Q_{(H)})_{(L)}\neq\varnothing$ if and only if $(H)\succ (L)$.
\end{enumerate}
\end{lm}

\begin{pr}\label{propofirst}
Suppose $G$ acts properly on the manifold $Q$. Then we have:

\begin{itemize}

\item[(i)] For $q\in Q_{(H)}$ such that $G_q=H$ and $(L)\in
I_{S^*Q}$,
$$
(J^{-1}(0))_{(L)}\cap S^*_qQ\neq\varnothing\iff (L)\in
I_Q\,\text{and}\,\left((H)\in I^*_Q\,
\text{or}\,(L)\prec(H)\right) ;
$$

\item[(ii)]  $(L)\in I_{J^{-1}(0)}\iff$ $(L)\in I^*_Q$ and hence
$\mathcal C_0^{(L)}\neq\varnothing\iff (L)\in I^*_Q\iff \dim
Q^{(L)}\ge 1$;

\item[(iii)]  The cosphere bundle projection $k$ restricts to the
$G$-equivariant continuous surjection
$k_{(L)}:(J^{-1}(0))_{(L)}\rightarrow \ov{Q_{(L)}}$ which is also
an open map;

 \item[(iv)]  For a fixed orbit type $(L)$ in the
zero momentum level set of the lifted $G $-action to $S^{*}Q$ the
corresponding orbit type submanifold admits the following
$G$-invariant partition:
\begin{equation}
\label{unu} (J^{-1}(0))_{(L)}=
J^{-1}_{(L)}(0)\coprod_{(H)\succ(L)}\pi_+
\left(J^{-1}_{ct(H)}(0)\times
\left(N^{*}Q_{(H)}\right)_{(L)}\right),
\end{equation}
where $(H)\in I_Q$;

\item[(v)]  For every $(H)\succ(L)$ with $(L)\in I^*_Q$ and
$(H)\in I_Q$ the restrictions
\begin{equation*}
\widetilde{t}_{(L)}:=k_{(L)}|_{J^{-1}_{(L)}(0)}\quad
\textrm{and}\quad \widetilde{t}_{(H)\succ
(L)}:=k_{(L)}|_{\pi_+\left({J^{-1}_{ct(H)}(0)}\times
(N^{*}Q_{(H)})_{(L)}\right)}
\end{equation*}
are $G$-equivariant smooth surjective submersions onto $Q_{(L)}$
and $Q_{(H)}$ respectively. The mappings $J_{ct(H)}$ and $J_{(H)}$
denote the momentum maps of the restricted actions of $G$ to
$T^{*}Q_{(H)}$ and $S^{*}Q_{(H)}$ respectively (which are the same
as the canonical momentum maps for the restricted $G$-action on
$Q_{(H)})$.
\end{itemize}
\end{pr}

\begin{proof}
To prove (i), let $(L)\in I_{S^*Q}$ and $q\in Q_{(H)}$ with
$G_q=H$. Then
\begin{equation}
(J^{-1}_{ct}(0))_{(L)}\cap T^*_qQ=(S_q^{H})^*\oplus
(N_q^*Q_{(H)})_{(L)},
\end{equation}
where $S_q^H$ is the linear slice for the $G$-action on
$Q_{(H)}$(see section 3 of \cite{prs}). Since
$(J^{-1}(0))_{(L)}\cap S^*_qQ=\varnothing\iff
(J^{-1}_{ct}(0))_{(L)}\cap T^*_qQ=\{0\}$, then
$(J^{-1}(0))_{(L)}\cap S^*_qQ=\varnothing$ only when $(S_q^{H})^*$
and $(N_q^*Q_{(H)})_{(L)}$ are simultaneously zero. This amounts
to $(L)=(H)\in I_Q\backslash I_Q^*$, (see Lemma
\ref{lemmacotlift}) from where the result follows.

(ii) is a forward consequence of (i). The rest of this statement
and the $G$-equivariant continuous surjectivity of $k_{(L)}$ are
 direct consequences of the fact that $I_Q=I_{J_{ct}^{-1}(0)}$.
 To prove the openness of $k_{(L)}$ it suffices to observe that
 for any open subset $\mathcal U$ of $(J^{-1}(0))_{(L)}$,
 $k_{(L)}(\mathcal U)= \tau_{(L)}(\pi^{-1}(\mathcal U ))$, where
 $\tau_{(L)}:(J_{ct}^{-1}(0))_{(L)}\rightarrow\ov{Q_{(L)}}$ is the
 open canonical cotangent projection map.

Applying \eqref{descomp} and the fact that $(N^{*}Q_{(H)})_{(L)}$
does not contain the zero section when $(H) \neq (L)$
  we have
\begin{equation*}
(J^{-1}_{ct}(0)\setminus \{0_{T^{*}Q}\})_{(L)} = \left(
J^{-1}_{ct(L)}(0)\right) \setminus
\{0_{T^{*}Q_{(L)}}\}\coprod_{(H)\succ(L)}
\left[J^{-1}_{ct(H)}(0)\times
\left(N^{*}Q_{(H)}\right)_{(L)}\right].
\end{equation*}
Hence, applying $\pi_+$ to this relation, we get
$$
(J^{-1}(0))_{(L)}=J^{-1}_{(L)}(0)\coprod_{(H)\succ(L)}
\pi_+\left(J^{-1}_{ct(H)}(0)\times
\left((N^{*}Q_{(H)}\right)_{(L)}\right)
$$
which proves statement (iv).

As for the proof of (v), it is enough to notice that
\[
J^{-1}_{(L)}(0) \quad \text{and} \quad
\pi_+\left(J^{-1}_{ct(H)}(0)\times(N^{*}Q_{(H)})_{(L)}\right)
\]
are bundles over $Q_{(L)}$ and $Q_{(H)}$ respectively.
\end{proof}

\begin{re}
Notice that for the description of orbit types in $J^{-1}(0)$, we
need not only $I_Q$, but also the lattice $I^*_Q$ since each
$J^{-1}(0)_{(L)}$ is written as a union with index $(H)$
in $I_Q$, but $(L)$ belongs to $I_Q^*$.
\end{re}

\section{Topology and contact geometry of $\mathcal{C}_0$}

\subsection{The secondary decomposition of $\mathcal{C}_0^{(L)}$}

Define the fiber bundles:
$$
s_{(H) \succ (L)}:=J^{-1}_{ct(H)}(0)\times (N^{*}Q_{(H)})_{(L)}
\rightarrow Q_{(H)}
$$
\[
s_{(L)}:=J^{-1}_{ct(L)}(0)\rightarrow Q_{(L)}.
\]
Taking into account that $\pi_+(s_{(H) \succ (L)})$ are
 $G$-invariant
pieces of the partition \eqref{unu} of $(J^{-1}(0))_{(L)}$ and
that the actions of $G$ and $\RR_+$ commute, we can define:

  $$
CS_{(H) \succ (L)}:=\frac{ \pi_+(s_{(H) \succ (L)})}G
$$
\[
CC_{(L)}:=\frac{J_{(L)}^{-1}(0)}G=\frac{\pi_+\left(s_{(L)}
\setminus \{0_{T^*Q_{(L)}}\}\right)}{G}\simeq
S^{*}\left(\frac{Q_{(L)}}G\right).
\]
Notice that for the above equivalence we have applied Theorem
\ref{uf} and that each contact stratum admits the following
partition, which is the quotient of \eqref{unu}:

\begin{align}
 \label{secondary partitionc}
\mathcal
C_0^{(L)}&=\frac{\left(J^{-1}(0)\right)_{(L)}}G
=CC_{(L)}\coprod_{(H)\succ(L)}CS_{(H)
\succ (L)}  \\
&\simeq S^{*}\left(\frac
{Q_{(L)}}G\right)\coprod_{(H)\succ(L)}CS_{(H) \succ (L)}.
\nonumber
\end{align}

\begin{re}\label{remarksecondarymaps}
In the notations of the previous section, the maps $k_{(L)}$,
$\widetilde{t}_{(L)}$, and $\widetilde{t}_{(H)\succ (L)}$ descend
to
\[
k^{(L)}:\mathcal{C}_0^{(L)}\rightarrow\ov{Q^{(L)}}
 \text{,}\quad
 \widetilde{t}^{(L)}:CC_{(L)}\rightarrow Q^{(L)}
  \text{, and} \quad
\widetilde{t}^{(H)\succ (L)}: CS_{(H)\succ (L)}\rightarrow Q^{(H)};
\]
$k^{(L)}$ is an open continuous surjection and the other two are smooth surjective
submersions.
\end{re}

\begin{te}\label{theoremsecondary}
With the above notations, we obtain the following:
\begin{itemize}
\item[\rm(i)] $\overline{Q^{(L)}}$ is a stratified space with
strata $Q^{(H)}$, for all $(L)\preccurlyeq(H)$ and with frontier
condition given by
$$
Q^{(K)}\cap\ov{Q^{(H)}}\neq\varnothing\iff (H)\preccurlyeq(K).
$$ Moreover, $Q^{(L)}$ is open and dense in
$\ov{Q^{(L)}}$. \item[\rm(ii)] For every $(L)\in I^*_Q$ and
$(H)\in I_Q$, the partition \eqref{secondary partitionc} is a
stratification of the corresponding contact stratum $\mathcal
C_0^{(L)}$, called the \textit{secondary stratification}. The
frontier conditions are given by:
\[
CS_{(H)\succ(L)}\subset \partial CC_{(L)}\quad \text{for all} \quad
 (H)\succ(L);
\]
\[
CS_{(H')\succ(L)}\subset \partial CS_{(H)\succ(L)} \iff
 (H')\succ(H)\succ(L).
\]
Moreover, the piece $CC_{(L)}$ is diffeomorphic to
$S^*Q^{(L)}$, is open and dense in $\mathcal C_0^{(L)}$, and the
map $k^{(L)}$ is a surjective submersion of stratified spaces.
\end{itemize}
\end{te}
\begin{proof}
Since the $G$-action is proper, the orbit type decomposition of
$Q$ induces a stratification of $Q/G$ and the first part of the
theorem follows immediately considering the relative topology of
$\overline{Q^{(L)}}$ in $Q/G$. Also, \eqref{secondary partitionc}
is a locally finite partition and its pieces are obviously
submanifolds of $\mathcal C_0^{(L)}$. As $k^{(L)}$
 is a continuous map and $(k^{(L)})^{-1}(Q^{(L)})=CC_{(L)}$, it
 follows that $CC_{(L)}$ is open in $\mathcal C_0^{(L)}$. In order to prove
 the density, let $x\in\mathcal C_0^{(L)}$ and $\mathcal U$ be any
 open neighborhood of $x$. Hence, $\mathcal V= k^{(L)}(\mathcal
 U)$ is an open subset of $\ov {Q^{(L)}}$ and, since $Q^{(L)}$
 is dense in $\ov {Q^{(L)}}$, there is at least one element
 $y\in \mathcal V \cap Q^{(L)}$. Notice that
 $(k^{(L)})^{-1}(y)=(\widetilde{t}^{(L)})^{-1}(y)\subset CC_{(L)}$ and that there is at least an element
in $(\widetilde{t}^{(L)})^{-1}(y)$ which is in $\mathcal U$. This
means that $\mathcal U\cap CC_{(L)}\neq \varnothing$ which proves
the density of $CC_{(L)}$.

Using the density of $CC_{(L)}$, the first frontier condition for
the secondary stratification becomes obvious. For the second one,
consider in $\mathcal C_0^{(L)}$ an arbitrary open neighborhood
$\mathcal U$ of a point $x\in CS_{(H')\succ(L)}$. By the openness
property of $k^{(L)}$, we obtain that $\mathcal O=k^{(L)}(\mathcal
U)$ is an open neighborhood of $k^{(L)}(x)$ in $\ov {Q^{(L)}}$.
Applying (i), we have that $\mathcal O\cap
Q^{(H)}\neq\varnothing\iff(H')\succ(H)\succ(L)$. Furthermore, the
surjectivity of $\widetilde{t}^{(H)\succ (L)}$ implies
 $(\widetilde{t}^{(H)\succ (L)})^{-1}(z)\cap\mathcal
U\neq\varnothing$ for any $z\in\mathcal O\cap Q^{(H)}$, proving
that $CS_{(H')\succ(L)}\subset \partial CS_{(H)\succ(L)} \iff
 (H')\succ(H)\succ(L)$.

 As $k^{(L)}$ restricted to each piece of the secondary
decomposition is surjective, Remark \ref{remarksecondarymaps}
immediately implies that this map is a stratified surjective
submersion.
\end{proof}

We will refer to the strata of the form $CS_{(H)\succ(L)}$ as
\textit{contact seams} due to their stitching role that will be
explained later in Remark \ref{role}.

This theorem completes the topological description of each contact
stratum $\mathcal{C}_0^{(L)}$ in terms of its secondary
stratification. We shall now begin the investigation of
geometrical aspects, namely to what extent the strata of this
secondary stratification admit canonical contact structures in the
sense that the $1$-forms generating them are induced by some
cosphere bundle structures compatible with the reduced contact
form on the contact stratum. Thus, denote by
\[
\widetilde\Psi^{(H)}:CC_{(H)}\rightarrow \
(S^*Q^{(H)},\Theta^{(H)}_{\Sigma}\ )
\]
the bundle morphism given by Theorem \ref{uf}, where
$\Theta^{(H)}_{\Sigma}$ is a contact form on the cosphere bundle
of $Q^{(H)}$. Observe that the restricted projection onto the
first factor
\[
\mathrm{p}_{1(H)\succ (L)}:\left( J_{ct(H)}^{-1}(0)\setminus
\{0_{T^*Q_{(H)}}\}\right)\times (N^*Q_{(H)})_{(L)}\rightarrow
J_{ct(H)}^{-1}(0)\setminus \{0_{T^*Q_{(H)}}\}
\]
is $\RR_+$ and $G$-equivariant so it descends to
the surjective submersion
\[
\widetilde{\mathrm{p}}_1^{(H)\succ (L)}:CS^{\circ}_{(H)\succ(L)}
\rightarrow CC_{(H)},
\]
where
\[
CS^{\circ}_{(H)\succ(L)}:=\frac{\pi_+\left(
J_{ct(H)}^{-1}(0)\setminus \{0_{T^*Q_{(H)}}\}\times
(N^*Q_{(H)})_{(L)}\right)}{G}
\]
is an open and dense submanifold of  the contact seam
$CS_{(H)\succ(L)}$. Then, for any pair $(H)\succ(L)$, we have the
following bundle map covering the identity on $Q^{(H)}$
$$
\widetilde{\Psi}^{(H)\succ(L)}:=\widetilde{\Psi}^{(H)}\circ\widetilde{\mathrm{p}}_1^{(H)\succ
(L)}:CS^{\circ}_{(H)\succ(L)}\rightarrow S^*Q^{(H)}
$$
which is also a surjective submersion. We are now able to endow
each cosphere-like stratum $CC_{(H)}$ and each
$CS^{\circ}_{(H)\succ(L)}$ with $1$-forms given by:
\begin{equation}
\label{defeta}
\left(CC_{(H)},\eta_{(H)}:=(\widetilde{\Psi}^{(H)})^*\Th_\Sigma^{(H)}\right),\quad\text{and}\end{equation}
\begin{equation}\label{defeta2}
\left(CS^{\circ}_{(H)\succ(L)},\eta_{(H)\succ(L)}:=(\widetilde{\Psi}^
{(H)\succ(L)})^*\Th_\Sigma^{(H)}\right).
\end{equation}
It is impossible to induce in this way a $1$-form on the whole
piece $CS_{(H)\succ(L)}$ and hence we are forced to restrict
ourselves, for the time being, to $CS^{\circ}_{(H)\succ(L)}$. However,
we will show later how to extend this form to the whole
$CS_{(H)\succ(L)}$.

Theorem \ref{WL1} gives the existence of an abstractly defined
contact structure on each contact piece $C_0^{(L)}$ generated by a
$1$-form ${\theta_{\sigma}}_0^{(L)}$. One of the aims of this
section is to investigate the compatibility of the previously
defined forms $\eta_{(H)}$ and $\eta_{(H)\succ (L)}$ with the
reduced contact form ${\theta_{\sigma}}_0^{(L)}$ and to describe
as much as possible this abstract contact structure.

\begin{te}\label{theoremsecondaryb}
The strata $CC_{(L)}$ and $CS^\circ_{(H)\succ(L)}$ within the
contact stratum $\mathcal{C}_0^{(L)}$ satisfy the following
properties:
\begin {enumerate}
\item [(i)]  $(CC_{(L)},\eta_{(L)})$ is an open dense contact
submanifold of the contact stratum $\mathcal{C}_0^{(L)}$
contactomorphic to $(S^*(Q^{(L)}),\Theta_\Sigma^{(L)})$. \item
[(ii)]Using the above notations, the conformal classes of
$\eta_{(L)}$ and $\eta_{(H)\succ(L)}$ admit smooth extensions to
$\mathcal{C}_0^{(L)}$ equivalent to ${\theta_{\sigma}}_0^{(L)}$,
namely
$$
{\theta_{\sigma}}_0^{(L)}\rrestr{{CC_{(L)}}}\simeq\eta_{(L)}\quad\text{and}\quad
{\theta_{\sigma}}_0^{(L)}\rrestr{{CS^\circ_{(H)\succ(L)}}}\simeq\eta_{(H)\succ(L)}.
$$
The extension of $\eta_{(L)}$ is unique.

\item[(iii)] The conformal class of $\eta_{(H)\succ(L)}$ can be
smoothly and uniquely extended to the whole stratum
$CS_{(H)\succ(L)}$. If $(H)\in I_Q^*$ then $CS_{(H)\succ(L)}$ is a
coisotropic submanifold of the contact stratum $\mathcal{C}_0^{(L)}$. When  $(H)\in I_Q\setminus I_Q^*$ then $CS_{(H)\succ(L)}$ is a
Legendrian submanifold of the contact stratum $\mathcal{C}_0^{(L)}$.
\end{enumerate}
\end{te}
\begin{proof}
(i) is a simple consequence of Theorem \ref{uf}.

For (ii), let $(L)$ and $(H)$ be two fixed elements of $I_Q^*$ and
$I_Q$ respectively and
$i_0^{(H)\succ(L)}:CS^\circ_{(H)\succ(L)}\rightarrow
\mathcal{C}_0^{(L)}$ the inclusion map. By definition,
$$\begin{array}{lll}
{\theta_{\sigma}}_0^{(L)}\rrestr{{CS^{\circ}_{(H)}\succ(L)}}\simeq\eta_{(H)\succ(L)}
 & \iff & \exists f>0 \quad\text{in}\quad C^\infty(CS_{(H)\succ (L)}) \quad\text{such
 that}\\

{\theta_{\sigma}}_0^{(L)}\rrestr{{CS^{\circ}_{(H)\succ(L)}}}=f\,\eta_{(H)\succ(L)}

& \iff &
 (i_\circ^{(H)\succ(L)}\ )^*{\theta_{\sigma}}_0^{(L)}\simeq\
(\widetilde{\Psi}^ {(H)\succ(L)})^*\Th_\Sigma^{(H)}.\end{array}
$$
 To simplify the reading of the proof, consider the two
 figures \ref{diag1} and \ref{diag2} where $\pi_G^{(H)\succ(L)}$ and $\bar{\pi}^{(H)}_G$ denote the
canonical $G$-projections and all the horizontal arrows in the
first and second diagram are injections and projections
respectively.

\begin{figure}

\small$$\small\xymatrix{ J_{ct(H)}^{-1}(0)\setminus
\{0_{T^*Q_{(H)}}\}\times (N^*Q_{(H)})_{(L)}
\ar@{^{(}->}[r]^>{\hspace{-6mm}i_{(H)\succ(L)}} \ar[d]_{\pi_+} &
\, (J_{ct}^{-1}(0)\setminus\{0_{T^*Q}\})_{(L)} \ar[d]^{\pi_+}
\ar@{^{(}->}[r]^>{\hspace{-7mm}i_{(L)}} & T^*Q\backslash
\{0_{T^*Q}\}\hspace{95mm}\ar[d]^{\pi_+}\hspace{-90mm} &
\\
 \small\pi_+\left( J_{ct(H)}^{-1}(0)\setminus \{0_{T^*Q_{(H)}}\}\times
(N^*Q_{(H)})_{(L)}\right)\ar@{^{(}->}[r]^>{\hspace{-12mm}\widetilde
i_{(H)\succ(L)}}\ar[d]_{\pi_G^{(H)\succ(L)}} &
(J^{-1}(0))_{(L)}\ar@{^{(}->}[r]^{\widetilde
i_{(L)}}\ar[d]^{\pi_G^{(L)}}& S^*Q\\
CS^{\circ}_{(H)\succ(L)}\ar@{^{(}->}[r]^{i_0^{(H)\succ(L)}}\ar[d]_{\widetilde{\Psi}^
{(H)\succ(L)}} & \mathcal{C}_0^{(L)}\\
S^*Q^{(H)}}$$
 \caption{Diagram defining $\eta_{(H)}$}
 \label{diag1}\end{figure}
\begin{figure}
\small
$$ \xymatrix{\\\left( J_{ct(H)}^{-1}(0)\setminus
\{0_{T^*Q_{(H)}}\}\times
(N^*Q_{(H)})_{(L)}\right)\ar[r]^{\hspace{12mm}p_{1(H)\succ(L)}}\ar[d]_{\pi_+}
&
J_{ct(H)}^{-1}(0)\setminus \{0_{T^*Q_{(H)}}\}\ar[d]^{\pi_+} &\\
\pi_+\left( J_{ct(H)}^{-1}(0)\setminus \{0_{T^*Q_{(H)}}\}\times
(N^*Q_{(H)})_{(L)}\right)\ar[r]^>{\hspace{-15mm}\widetilde{p}_{1(H)\succ(L)}}\ar[d]_{\pi_G^{(H)\succ(L)}}
&
J_{(H)}^{-1}(0)\ar[d]^{\bar{\pi}_G^{(H)}}\ar[r]^{l_{(H)}} & S^*Q_{(H)} \\
CS^\circ_{(H)\succ(L)}\ar[r]^{\widetilde{p_1}^{(H)\succ(L)}} &
CC_{(H)}\ar[r]^{\widetilde{\psi}^{(H)}} & S^*Q^{(H)} }
$$
\caption{Diagram defining $\eta_{(H) \succ (L)}$}\label{diag2}\end{figure}

As $\pi_G^{(H)\succ(L)}\circ\pi_+$ is a submersion, it suffices to
prove that
\begin{equation}
\label{ecu}
 \ ( i_0^{(H)\succ(L)}\circ\pi_G^{(H)\succ(L)}\circ\pi_+
\ )^*{\theta_{\sigma}}_0^{(L)}\simeq\ (\widetilde{\Psi}^
{(H)\succ(L)}\circ\pi_G^{(H)\succ(L)}\circ\pi_+\
)^*\Th_\Sigma^{(H)}.
\end{equation}
Observe that $i_0^{(H)\succ(L)}\circ\pi_G^{(H)\succ(L)}\circ\pi_+=
\pi_G^{(L)}\circ\pi_+\circ i_{(H)\succ(L)}$ and the first term of
\eqref{ecu} becomes
\begin{align*}
\ (\pi_G^{(L)}\circ\pi_+\circ i_{(H)\succ(L)}\
)^*{\theta_{\sigma}}_0^{(L)}=i_{(H)\succ(L)}^*\circ\pi_+^*\
((\pi_G^{(L)})^*{\theta_{\sigma}}_0^{(L)}\ )\\=
i_{(H)\succ(L)}^*\circ\pi_+^*\ (\widetilde i_{(L)}^*\theta_\sigma\
)=\ (\pi_+\circ j_{(H)}\circ\Phi\ )^*\theta_\sigma,
\end{align*}
where in the last line we have used Theorem \ref{WL1} together
with the equality $\pi_+\circ j_{(H)}\circ\Phi=\widetilde
i_{(L)}\circ\pi_+\circ i_{(H)\succ(L)}$, with $j_{(H)}$ and $\Phi$
inclusions defined by:
$$
\Phi:\left( J_{ct(H)}^{-1}(0)\setminus \{0_{T^*Q_{(H)}}\}\times
(N^*Q_{(H)})_{(L)}\right)\hookrightarrow T^*Q|_{Q_{(H)}}\setminus
\{0_{T^*Q_{(H)}}\} $$
 and
$$
 j_{(H)}:T^*Q|_{Q_{(H)}}\setminus
\{0_{T^*Q_{(H)}}\}\hookrightarrow T^*Q \setminus \{0_{T^*Q}\}.
 $$
Using this time $\bar{\pi}_G^{(H)}\circ\pi_+\circ
p_{1(H)\succ(L)}= \widetilde
p_1^{(H)\succ(L)}\circ\pi_G^{(H)\succ(L)}\circ\pi_+$, we can write
the second term of \eqref{ecu} as:
\begin{align*}
\ (\widetilde{\Psi}^ {(H)}\circ\widetilde p_1^{(H)\succ(L)}\circ
\pi_G^{(H)\succ(L)}\circ\pi_+\ )^*\Th_\Sigma^{(H)}=\
(\widetilde{\Psi}^ {(H)}\circ\bar{\pi}_G^{(H)}\circ\pi_+\circ
p_{1(H)\succ(L)}\ )^*\Th_\Sigma^{(H)}\\=\ (\pi_+\circ
p_{1(H)\succ(L)}\ )^*\ (\widetilde{\Psi}^
{(H)}\circ\bar{\pi}_G^{(H)}\ )^*\Th_\Sigma^{(H)}\simeq\
(\pi_+\circ p_{1(H)\succ(L)}\ )^*l_{(H)}^*\theta_{(H)\Sigma},
\end{align*}
where  $\theta_{(H)\Sigma}$ is a contact form on $S^*Q_{(H)}$. Let
$p_{(H)}:T^*Q|_{Q_{(H)}}\setminus \{0_{T^*Q_{(H)}}\}\rightarrow
T^*Q_{(H)}$ be the projection map. Since $l_{(H)}\circ\pi_+\circ
p_{1(H)\succ(L)}=\pi_+\circ p_{(H)}\circ\Phi$, the second term is
in the same conformal class as $\Phi^*
p_{(H)}^*\pi_+^*\theta_{(H)\Sigma}$ and, hence, equation
\eqref{ecu} is equivalent to
\begin{align*}
\Phi^*p_{(H)}^*\pi_+^*\theta_{(H)\Sigma}\simeq
\Phi^*j_{(H)}^*\pi_+^*\theta_\sigma\iff
\Phi^*p_{(H)}^*(\Sigma\circ\pi_+)^*\theta_{(H)}\simeq
\Phi^*j_{(H)}^*(\sigma\circ\pi_+)^*\theta\\
\iff \Phi^*p_{(H)}^*f_\Sigma\theta_{(H)}\simeq
\Phi^*j_{(H)}^*f_\sigma\theta ,
\quad\quad\quad\quad\quad\quad\quad
\end{align*}
where $\theta$ and $\theta_{(H)}$ are the canonical one-forms on
$T^*Q$ and $T^*Q_{(H)}$ respectively, and $\sigma$, $\Sigma$ are
sections in the associated cosphere bundles. But
$p_{(H)}^*\theta_{(H)}=j_{(H)}^*\theta$, as it can be easily seen
in local coordinates, which proves \eqref{ecu}.

As for the extension of the conformal class of $\eta_{(L)}$, an
analogous proof can be developed just by considering the limit
case $(H)=(L)$, when $CS^0_{(H)\succ (L)}$ degenerates in
$CC_{(L)}$. In order to prove the uniqueness of this extension,
let us consider a point $x\in\mathcal{C}_0^{(L)}$ and one tangent
vector $v_x\in T_x\mathcal{C}_0^{(L)}$. As $CC_{(L)}$ is open and
dense in $\mathcal{C}_0^{(L)}$, there is  a sequence of points
$x_k\in CC_{(L)}$ and one of vectors $v_{x_k}\in
T_{x_k}CC_{(L)}\simeq T_{x_k}\mathcal{C}_0^{(L)}$ such that
\[
\lim_{k\to\infty}x_k =x\text{,}\quad\lim_{k\to\infty}v_{x_k}=v_x.
\]
From the above arguments and using the continuity of
${\theta_{\sigma}}_0^{(L)}$, we have that
\[
\lim_{k\to\infty}\frac{\eta_{(L)}(x_k)(v_{x_k})}{g(x_k)}=
\lim_{k\to\infty}{\theta_{\sigma}}_0^{(L)}(x_k)(v_{x_k})=
{\theta_{\sigma}}_0^{(L)}(x)(v_x),
\]
with $g\in C^\infty(CC_{(L)})$ a positive function such that
$\eta_{(L)}=g\,{\theta_{\sigma}}_0^{(L)}\rrestr{{CC_{(L)}}}$. We
have thus proved that the class of ${\theta_{\sigma}}_0^{(L)}$ is
the unique smooth extension of the class of $\eta_{(L)}$ to
$\mathcal{C}_0^{(L)}$.

(iii) To extend the class of $\eta_{(H)\succ(L)}$ from
$CS^\circ_{(H)\succ(L)}$ to the whole piece $CS_{(H)\succ(L)}$, we
will apply the same type of arguments as before, using this time
that $CS^\circ_{(H)\succ(L)}$ is open and dense in
$CS_{(H)\succ(L)}$. Namely, for any point $x\in CS_{(H)\succ(L)}$
and any $v_x\in T_{x}CS_{(H)\succ(L)}$, there is  a sequence of
points $x_k\in CS^\circ_{(H)\succ(L)}$ and one of vectors
$v_{x_k}\in T_{x_k}CS^\circ_{(H)\succ(L)}\simeq
T_{x_k}CS_{(H)\succ(L)}$ such that
\[
\lim_{k\to\infty}x_k =x\text{,}\quad\lim_{k\to\infty}v_{x_k}=v_x.
\]
Observe that
 \[
\lim_{k \to\infty}\frac{\eta_{(H)\succ(L)}(x_k)(v_{x_k})}{f(x_k)}=
{\theta_{\sigma}}_0^{(L)}(x)(v_x)
\]
and notice that this extension is also unique and given by the
conformal class of
${\theta_{\sigma}}_0^{(L)}|_{CS_{(H)\succ(L)}}$.

 To check the coisotropy and Legendrian submanifold conditions, let $x\in
 CS^\circ_{(H)\succ(L)}$. A direct count of dimensions
 gives:
 $$
\operatorname{dim} \ker{\theta_{\sigma}}_0^{(L)}(x)
=\operatorname{dim}\mathcal{C}_0^{(L)}-1=
 2(\operatorname{dim}Q_{(L)}-\operatorname{dim}G +\operatorname{dim} L -1)$$
since $S^*Q^{(L)}$ is open in the corresponding contact stratum.
At this point we need the following intermediate result.
\begin{lm} The dimension of the tangent space to a contact seam is
\begin{equation}\label{dimensionseam} \operatorname{dim}T_xCS_{(H)\succ(L)}= \operatorname{dim}
Q_{(H)}+ \operatorname{dim} Q_{(L)}-2 \operatorname{dim} G+
\operatorname{dim} H+ \operatorname{dim} L-1. \end{equation}
\end{lm}
\begin{proof}
We want to compute $\operatorname{dim}T_xCS_{(H)\succ(L)}
=\operatorname{dim} CS_{(H)\succ(L)}$. For this, let $\pi(z)=k^0
(x)$ be the base point of $x$, where $z\in Q_{(H)}$ with $G_z=H$
and note that $\operatorname{dim}
CS_{(H)\succ(L)}=\operatorname{dim}(J_{ct}^{-1}(0)\cap
T_{z}^{\ast}Q)_{(L)}+\operatorname{dim}
Q_{(H)}-\operatorname{dim}G+\operatorname{dim}L-1$. Where the
class $(L)$ refers to the linear $H$-action on the vector space
$J_{ct}^{-1}(0)\cap T_{z}^{\ast}Q$. On the other hand, the inverse
of the Riemannian bundle isomorphism $TQ\rightarrow T^*Q$ maps
$(J_{ct}^{-1}(0)\cap T^{\ast}_zQ)_{(L)}$ $H$-equivariantly
isomorphically to $(S_z)_{(L)}$. Now, if $\psi$, $U$, and $U'$ are
like in the Tube Theorem \eqref{tube}, then $\psi$ restricts to a
diffeomorphism between $G\times_{H} \left((S_z)_{(L)}\cap
U\right)$ and $U\cap Q_{(L)}$. Since
$\operatorname{dim}G\times_{H}(S_{z})_{(L)}=\operatorname{dim}G+\operatorname{dim}(S_{z})_{(L)}-\operatorname{dim}
H$, we can compute
\begin{equation*}
\operatorname{dim}(S_z)_{(L)}=\operatorname{dim}
Q_{(L)}-\operatorname{dim} G+\operatorname{dim} H.
\end{equation*}
Finally we obtain $\operatorname{dim}T_xCS_{(H)\succ(L)}
=\operatorname{dim}Q_{(H)}+\operatorname{dim}Q_{(L)}-2\operatorname{dim}G
+\operatorname{dim}H +\operatorname{dim}L-1$.
\end{proof}
 Consequently, a simple dimension count gives
\begin{align*}
\operatorname{dim} T_xCS_{(H)\succ(L)} -
\frac{1}{2}\operatorname{dim} \ker{\theta_{\sigma}}_0^{(L)}(x)
&=\operatorname{dim} Q_{(H)}- \operatorname{dim} G+
\operatorname{dim}
H\\
&=\operatorname{dim}(S_z)_{(H)}\geq 0,
\end{align*}
where $z \in Q_{(H)}$ is the base point of $x $ and $S_z $ is the
associated linear slice.
 Suppose first that $(H)\in I_Q^*$ and so $\operatorname{dim} T_xCS_{(H)\succ(L)} -
\frac{1}{2}\operatorname{dim}
\ker{\theta_{\sigma}}_0^{(L)}(x)\gneqq 0$. This implies that
$CS^\circ_{(H)\succ(L)}$ and $CS_{(H)\succ(L)}$ can be neither
isotropic nor Legendrian submanifolds of $\mathcal{C}_0^{(L)}$ and that
 $T_xCS^\circ_{(H)\succ(L)}\varsubsetneqq\ker
{\theta_{\sigma}}_0^{(L)}(x)$ for any $x\in
CS^\circ_{(H)\succ(L)}$.

Now let
$$W_x:=T_xCS^\circ_{(H)\succ(L)}\cap \ker
{\theta_{\sigma}}_0^{(L)}(x) =T_xCS_{(H)\succ(L)}\cap \ker
{\theta_{\sigma}}_0^{(L)}(x)$$ and $$V_x:=\left\{v\in
T_xCS^\circ_{(H)\succ(L)}\setminus\ker
{\theta_{\sigma}}_0^{(L)}(x):v=v_0\oplus kR(x),
k\in\RR\,\text{,}\,v_0\in\ker
{\theta_{\sigma}}_0^{(L)}(x)\right\}.$$ One can easily check that
$V_x$ is a one dimensional vector space and that for any $x\in
CS^\circ_{(H)\succ(L)}$, we have $T_xCS^\circ_{(H)\succ(L)}=W_x\oplus V_x$.  As $\widetilde{\Psi}^ {(H)\succ(L)}$ is a surjective submersion
and
${\theta_{\sigma}}_0^{(L)}\rrestr{{CS^\circ_{(H)\succ(L)}}}\simeq\eta_{(H)\succ(L)}$,
it follows that $T_x \widetilde{\Psi}^
{(H)\succ(L)}(W_x)=\ker\Theta_\Sigma^{(H)}(y)$ and $T_x
\widetilde{\Psi}^ {(H)\succ(L)}(V_x)$= span $\{R_\Sigma(y)\}$,
where $y=\widetilde{\Psi}^ {(H)\succ(L)}(x)$ and $R_\Sigma(y)$ is
the Reeb vector field of $(S^*Q^{(H)}, \Theta_\Sigma^{(H)})$.
Therefore, we obtain
\begin{align*}
&\operatorname{rank} d\eta_{(H)\succ(L)}(x)|_{W_x} =
\dim W_x-\operatorname{dim} \ker d\eta_{(H)\succ(L)}(x)|_{W_x} \\
&\quad = \operatorname{dim}W_x-\operatorname{dim} \{v\in
W_x:d\Theta_\Sigma^{(H)}(y)(T_x\widetilde{\Psi}^
{(H)\succ(L)}v,T_x\widetilde{\Psi}^ {(H)\succ(L)}w)=0, \forall
w\in W_x \} \\
&\quad = \operatorname{dim}W_x-\operatorname{dim}\operatorname{ker}T_x\widetilde{\Psi}^
{(H)\succ(L)}\rrestr{W_x} =
\operatorname{dim}S^*Q^{(H)}-1.\qquad\qquad\qquad
\end{align*}
This shows that $\operatorname{rank}
d\eta_{(H)\succ(L)}(x)\rrestr{W} = 2 \operatorname{dim}
W-(\operatorname{dim} \mathcal{C}_0^{(L)}-1)$ proving that
$CS^\circ_{(H)\succ(L)}$ is a coisotropic submanifold. Since
$CS^\circ_{(H)\succ(L)}$ is dense in $CS_{(H)\succ(L)}$, by an
extension argument similar to the one used before, we have that
$CS_{(H)\succ(L)}$ is also a coisotropic submanifold of the
corresponding contact stratum.

If $(H)\in I_Q\setminus I^*_Q$, then $\operatorname{dim}
T_xCS_{(H)\succ(L)}=\frac{1}{2}\operatorname{dim}
\ker{\theta_{\sigma}}_0^{(L)}(x)$ and by the definition
\eqref{defeta2} $\eta_{(H)\succ(L)}=0$ since $S^*Q^{(H)}$ is the
trivial bundle, proving thus that the $CS_{(H)\succ(L)}$ is a
Legendrian submanifold of $\mathcal{C}_0^{(L)}$.
\end{proof}

\begin{re}
Note that the contact seams $CS_{(H)\succ(L)}$ can never be
contact submanifolds of $\mathcal{C}_0^{(L)}$.
\end{re}

\subsection{The C-L stratification of $\mathcal{C}_0$}
In this subsection we prove the existence of a new stratification
of the contact reduced space $\mathcal{C}_0$, different from the
contact stratification in Theorem \ref{WL1}. The existence of this
new stratification, that we call the C-L stratification since its
strata are coisotropic or Legendrian submanifolds of the
corresponding contact stratum, is due to the bundle structure of
the contact manifold that we start with. We will see that the C-L
stratification is strictly finer than the contact one, if the base
manifold $Q$ has more than one orbit type. In principle, this is
not an advantage since the contact stratification partitions the
singular contact quotient in less and larger smooth components.
However, if we take into account the bundle structure of the
problem we can see why this new stratification is more
appropriate.

The most important feature of regular cosphere bundle reduction,
Theorem \ref{olt1}, is that if we start with the cosphere bundle
of a manifold $Q$, we end up again with a cosphere bundle, this time
over $Q/G$. Furthermore, the reduced contact structure on
$S^*(Q/G)$ equals the canonical cosphere contact structure. In the
singular setting however, the lack of smoothness of the
quotient spaces involved forces us to choose another definition of
fibration. The most natural one when working with decomposed or
stratified spaces is the following: if $A$ and $B$ are decomposed
spaces together with a continuous surjection $f:A\rightarrow B$,
we say that $f:A\rightarrow B$ defines a \textit{stratified
bundle} over $B$ if $f$ is a morphism of decomposed spaces. In our
case, there is a natural projection $k^0:\mathcal{C}_0\rightarrow
Q/G$ induced from the cosphere bundle projection
$k:S^*Q\rightarrow Q$. If we consider the natural orbit type
stratification of $Q/G$ and the contact one of $\mathcal{C}_0$,
then the projection does not define a stratified bundle over $Q/G$
since the image of a contact stratum $\mathcal{C}_0^{(L)}$ under
the projection is $\overline{Q^{(L)}}$ which includes several
orbit type strata of $Q/G$. We will prove that the choice of the
coisotropic stratification for the contact quotient
$\mathcal{C}_0$ solves this problem.

Consider the partition of $\mathcal{C}_0$ obtained by putting
together all the secondary strata found in every contact stratum:
\begin{equation}\label{coisotropicpartition}\mathcal{C}_0=\coprod_{(L)}CC_{(L)}\coprod_{(K')\succ(K)}
CS_{(K')\succ (K)}
\end{equation}
for every pair of classes $(L),(K)\in I^*_Q$ and every $(K')\in I_Q$.

\begin{te}\label{coisotropicstratification}
The partition \eqref{coisotropicpartition} is a decomposition of
$\mathcal{C}_0$ inducing a stratification,  called the
C-L stratification, that satisfies the following properties:
\begin{itemize}
\item[(1)] If $Q/G$ is connected and $(L_0)$ is the principal
orbit type in $Q$, then $CC_{(L_0)}$ is open and dense in
$\mathcal{C}_0$.

\item[(2)] $k^0:\mathcal{C}_0\rightarrow Q/G$ is a stratified
bundle with respect to the C-L stratification of $\mathcal{C}_0$
and the orbit type stratification of $Q/G$.

\item[(3)] If $I_Q$ consists of more than one class, the C-L
stratification is strictly finer than the contact one, and they
are identical otherwise.

\item[(4)] The frontier conditions for the C-L stratification of
$\mathcal{C}_0$ are:
\begin{itemize}
\item[(i)] $CC_{(K)}\subset \partial CC_{(H)} \Longleftrightarrow
(H) \prec (K)$

\item[(ii)] $CS_{(K)\succ (H)}\subset\partial CC_{(H)}
\Longleftrightarrow (H) \prec (K)$

\item[(iii)] $C_{(K)}\subset \partial CS_{(K)\succ (H)}
\Longleftrightarrow (H)\prec (K)$

\item[(iv)] $CS_{(K')\succ (H)}\subset \partial CS_{(K)\succ (H)}
\Longleftrightarrow (H)\prec (K)\prec (K')$

\item[(v)] $CS_{(K)\succ(H')}\subset\partial CS_{(K)\succ
(H)}\Longleftrightarrow (H)\prec (H')\prec (K).$
\end{itemize}
\end{itemize}
\end{te}

\begin{proof}
For (1), recall by Proposition \ref{propofirst} that $I_{J^
{-1}(0)}=I^*_Q$. The principal orbit type of the isotropy lattice
corresponds to an open and dense piece, so $(J^{-1}(0))_{(L_0)}$
is open and dense in $J^{-1}(0)$, since $(L_0)$ is by hypothesis
the principal orbit type in $I^*_Q$ (assuming that $\dim Q \neq
0$) and hence in $I_{J^ {-1}(0)}$. Consequently, as the orbit map
$J^{-1}(0)\rightarrow \mathcal{C}_0$ is continuous and open,
$\mathcal{C}_0^{(L_0)}$ is open and dense in $\mathcal{C}_0$. Now,
since $\mathcal{C}_0^{(L_0)}$ is equipped with the relative
topology with respect to $\mathcal{C}_0$ and $CC_{(L_0)}$ is open
and dense in it (Theorem \ref{theoremsecondary}), it follows that
$CC_{(L_0)}$ is also open and dense in $\mathcal{C}_0$.  For
(2), note that the restrictions of $k^0$ to $CC_{(L)}$ and
$CS_{(H)\succ (L)}$ coincide with the corresponding restrictions
of $k^{(L)}$, which, by Remark \ref{remarksecondarymaps}, are
smooth surjective submersions over $Q^{(L)}$ and $Q^{(H)}$
respectively for every $(L)\in I^*_Q$ and $(H)\in I_Q$. This shows
that these restrictions map each C-L stratum of $\mathcal{C}_0$ to
an orbit type stratum of $Q/G$. Therefore, $k^0$ is a morphism of
stratified spaces. To prove $(3)$, recall from Theorem \ref{uf},
that if $I_Q$ consists of a single orbit type $(H)$, then
$\mathcal{C}_0=\mathcal{C}_0^{(H)}=CC_{(H)}$ (assuming $\dim Q
\neq 0$) and its contact and C-L stratifications are both trivial
and identical. If there is more than one orbit type in the base,
the number of C-L strata is strictly higher than the number of
contact strata (which is equal to the number of orbit types of
$I^*_Q$). The identity map in $\mathcal{C}_0$ injects each C-L
stratum in the unique contact stratum to which it belongs and is
hence a morphism of stratified spaces. Therefore, the C-L
stratification is finer than the contact one. For $(4)$, relations
$(ii)$ and $(iv)$ follow from the frontier conditions of the
secondary stratum $\mathcal{C}_0^{(H)}$. To prove $(i)$, it
suffices to recall from the general theory of singular contact
reduction that
$\mathcal{C}_0^{(K)}\subset\partial\mathcal{C}_0^{(H)}$ if and
only if $(H)\prec (K)$. Using the density of any maximal secondary
stratum $CC_{(L)}$ in the corresponding contact piece
$\mathcal{C}_0^{(L)}$, $(i)$ follows. $(iii)$ is a consequence of
$(v)$ if one considers the limit case $CC_{(K)}=CS_{(K)\succ
(K)}$.

Finally, to prove $(v)$, choose a point $[x]\in CS_{(K)\succ
(H')}\subset \mathcal{C}_0$ and an open neighborhood $[x]\in
O\subset \mathcal{C}_0$. We shall show that $O\cap CS_{(K)\succ
(H)}\neq\varnothing$ if $(H)\prec (H')\prec (K)$. Let $x\in
J^{-1}(0)$ be a preimage of $[x]$. We can assume without loss of
generality that $G_x=H'$ and that the projection of $x$, i.e. the
point $z=k(x)\in Q$, satisfies $G_z=K$. Let $U$ be the only open
$G$-saturated set in $J^{-1}(0)$ such that $U/G=O$. Then,
identifying $S^*Q$ with the unit bundle in $T^*Q$ via a
$G$-invariant metric on $Q$, we have that $x$ is a unit covector
lying in the subset of the cotangent fiber at $z$ given by
$(S_z^K)^*\oplus (N_z^*Q_{(K)})_{(H')}$. By the general properties
of linear representations of compact groups on vector spaces and
the property (iii) of cotangent-lifted actions in Lemma
\ref{lemmacotlift}, it follows that $p_2(U\cap T_z^*Q)\cap
(N_z^*Q_{(K)})_{(H)}\neq\varnothing$ for every compact subgroup
$H$ of $K$ such that $H\prec H'$ and
$(N_z^*Q_{(K)})_{(H)}\neq\varnothing$, i.e., $(H)\in I_Q$. Here,
$p_2$ is the linear projection $(S_z^K)^*\oplus
N_z^*Q_{(K)}\rightarrow N_z^*Q_{(K)}$. From this, it follows that
if $x'\in p_2(U\cap T_z^*Q))\cap (N_z^*Q_{(K)})_{(H)}$, then
$[x']\in O\cap CS_{(K)\succ (H)}$.
\end{proof}

\begin{re}\label{role}
The previous result shows that, identifying a stratum $CC_{(H)}$
with $S^*Q_{(H)}$ as shown in Theorem \ref{theoremsecondaryb},
 the reduced space
$\mathcal{C}_0$ is almost everywhere a collection of cosphere
bundles, one for each orbit type stratum of positive dimension in
$Q/G$. These cosphere bundles satisfy the same frontier conditions
as their bases, i.e., $S^*Q^{(K)}\subset \partial S^*Q^{(H)}$ if
and only if $Q^{(K)}\subset \partial Q^{(H)}$ (condition (i)), but
in this case there is always a contact seam $CS_{(K)\succ (H)}$
between them, which ``glues together'' these two cosphere bundles,
as reflected in conditions (ii) and (iii).

\end{re}

\section{Singular actions on the base with regular lifts to the cosphere bundle}

In the following definition we introduce a class of actions which
may have singularities on $Q$ but that will be proven to yield
regular lifted actions on $S^*Q$. \begin{de}\label{almostsemifree}
An \textit{almost semifree action} of $G$ on $Q$ is a smooth
action such that a) it is free almost everywhere, b) the connected
components of every orbit of non-maximal dimension are isolated,
and c) for every non-trivial isotropy subgroup $H\in I_Q$ with Lie
algebra $\mathfrak{h}$, its induced adjoint representation on
$\left(\mathfrak{g} / \mathfrak{h}\right)\setminus \{0\}$ given by
$h\cdot\left[\xi\right]=\left[\Ad_{h}\xi\right]$ is free.
\end{de} Note that for any almost semifree action, the quotient
space $Q/G$ consists on an open and dense stratum $Q^{(e)}$,
except possibly for a set of isolated singular points. The next
proposition shows that the class of almost semifree actions is in
one-to-one correspondence with the class of free actions on
$S^*Q$.

 \begin{pr}
 Let $S^*Q$ be the cosphere bundle of $Q$ endowed with the lift of
a proper action of a Lie group $G$ on $Q$. This lifted action is
free if and only if the action on $Q$ is almost semifree.
 \end{pr}

 \begin{proof}
 Recall that, identifying with the help of a $G$-invariant Riemannian metric $S^*Q$ with the unit bundle $SQ\subset TQ$ and $TQ$ with $T^*Q$, $G$ acts
freely on $S^*Q$ if and only if its tangent-lifted action on $TQ$
is free on the unit bundle, and hence if it is free away from the
zero section (since by linearity the lifted action intertwines the
fiber rescaling by non-zero factors). Let $q\in Q$
 with stabilizer $G_{q}=H\neq \{ e\}$, $S\subset T_{q}Q$
 a linear slice for the $G$-action at $q$ and $v=\xi_Q (q)+s\in T_{q }Q \setminus \{0\}$.
 Note that all the admissible $\xi$'s differ
by an element of the Lie algebra of $H$. Then
$U=G\cdot\mathrm{exp}_{q}(S)$ is a
 $G$-invariant neighborhood of the orbit $G\cdot
 q=G\cdot\mathrm{exp}_{q}(0)$ and there is an $H$-isomorphism $f:T_q Q\rightarrow \mathfrak{g}/ \mathfrak{h}\times S$
 given by $f(\xi_Q(q)+s)=(\left[\xi\right],s)$, where the $H$-invariance is with respect to the linear action on $T_qQ$ and the diagonal
 action on $\mathfrak{g}/ \mathfrak{h}\times S$ given by $h\cdot(\left[\xi\right],s)=(\left[\operatorname{Ad}_h\xi\right],h\cdot s)$.
 Consequently, $G_v=H_v=H_s\cap H_{\left[\xi\right]}$.

 Suppose first that the lifted action of $G$ on $S^*Q$ is free.
 Then any point $q'\in U\setminus G\cdot q$ can be written as
 $q'=g\cdot\exp_q(s)$ for some $0\neq s\in S$ with $g\in G$ and
 $G_{q'}=gH_sg^{-1}=\{e\}$, since $G_s=H_s=\{e\}$ as assumed
 above. Hence the $G$-action on $Q$ is almost semifree.

For $v = \xi_Q(q_0) \in T_{q_0} Q \setminus \{0\}$ with $\xi \in
\mathfrak{g}$, $\xi \notin  \mathfrak{h}$ we obtain that
$G_v=\{e\}=H_v= H_{\left[\xi\right]}$, thus proving that the
induced adjoint representation on $\left(\mathfrak{g} /
\mathfrak{h}\right)\setminus \{0\}$ is free.

To prove the converse implication,  let $v\in
T_{q}Q\setminus\{0\}$ as before, with $v=s+\xi_Q(q_0)$, where
$s\in S$ and $\xi\in\mathfrak{g}$. If $s$ is different from zero,
multiplying it if necessary by a positive scalar smaller than one,
we can guarantee that $G\cdot \mathrm{exp}_{q}(s)\subset
U\setminus G\cdot q$. Shrinking $U$ if necessary, we can guarantee
that all of the points in $U\setminus G\cdot q_0$ have trivial
isotropy, since the orbits of non maximal dimension are isolated
by hypothesis. Using again the Tube Theorem, the isotropy groups
of these points are $gH_sg^{-1}=\{e\}$, for every $g\in G$, which
forces $H_s=\{e\}$ and hence $H_v=\{e\}$. In the case when $s=0$,
we have that $G_v=H_{\left[\xi\right]}=\{e\}$, thus completing the
proof.
\end{proof}

\begin{re}
To geometrically express the third condition in Definition
\ref{almostsemifree}, notice that every non trivial isotropy
subgroup $H=G_q\in I_Q$ acts freely on $\left(\mathfrak{g} /
\mathfrak{h}\right)\setminus \{0\}$ if and only if for any element
$h\in H$ the associated diffeomorphism of $Q$ maps bijectively
$\{\exp_q(t \xi) \cdot q \,:\, t \in \mathbb{R}\}$ to $\{\exp_q(t
\operatorname{Ad}_h \xi) \cdot q \,:\, t \in \mathbb{R}\}$ for
every $\xi\in\mathfrak{g}$ with $[\xi]\neq 0$ in
$\g/\mathfrak{h}$.
\end{re}

Notice that this is a major difference with the cotangent bundle
case, where the cotangent-lifted action is free if and only if the
base action is free as well. In the context of cosphere bundle
reduction the reason for the special interest in semifree actions
and in finding necessary and sufficient conditions for the
freeness of the lifted cosphere action is the following. Given a
cosphere bundle $\mathcal{C}=S^*Q$ with the lift of a proper
almost semifree action on $Q$, if we ignore the bundle structure
of the contact manifold $\mathcal{C}$ we are in the hypothesis of
regular contact reduction, since $G$ acts freely, properly, and by
strong contactomorphisms on $\mathcal{C}$. Therefore, the contact
reduced space $\mathcal{C}_0$ is a well defined smooth contact
manifold.

 On the other hand, since the action on $Q$ is not free in general, we
cannot apply the main result on regular cosphere bundle reduction
of \cite{dor} (see Theorem \ref{olt1}) because in that case the
quotient $Q/G$ will not be a smooth manifold. In  fact,
 one expects $\mathcal{C}_0$ to be a smooth
 reduced manifold fibrating continuously over the topological stratified
 space $Q/G$, but this bundle description cannot be achieved by only applying the
 scheme of regular cosphere bundle reduction. However, the results
 of the previous section will allow us to provide such a ``stratified bundle"
 picture of the contact quotient $\mathcal{C}_0$.
 Indeed, we have the following result.
 \begin{te}
 Let $G$ be a Lie group acting properly and almost semifreely on $Q$ and by lifts on the
 cosphere bundle $S^*Q$ with contact momentum map $J:S^*Q\rightarrow \mathfrak{g}^*$. Write the orbit
 type
 decomposition of $Q/G$ as
 $$Q/G=Q^{(e)}\coprod_{(H)\in I_Q\setminus I^*_Q}*^{(H)},$$ where $Q^{(e)}=Q_{(e)}/G$ is open and dense
 in
 $Q/G$ and   each $*^{(H)}$ with $(H)\in I_Q\setminus I^*_Q $ is an isolated point
 of some lower dimensional stratum $Q^{(H)}$ with $(H)\succ (e)$, lying in the boundary of $Q^{(e)}$. Then the
quotient $\mathcal{C}_0=J^{-1}(0)/G$ is a
 smooth manifold which can be decomposed as
 \begin{equation}\label{semifreequotient}
 \mathcal{C}_0\simeq S^*Q^{(e)}\coprod_{(H)\in I_Q\setminus
 I^*_Q}CS^{(H)}
 \end{equation}
 where each $CS^{(H)}$ is a trivial
 bundle over $*^{(H)}$ and a connected submanifold of
 $\mathcal{C}_0$ lying in the boundary of $S^*Q^{(e)}$. Moreover, the manifolds $CS^{(H)}$ are Legendrian submanifolds of
$\mathcal{C}_0$ in one-to-one correspondence
 with the singular orbits of the $G$-action on $Q$ and
 have dimension
 $\mathrm{dim\,}Q-\,\mathrm{dim\,}G-1$.
 \end{te}

\begin{proof}
Since $J^{-1}(0)$ consists of a single orbit type $(e)$, due to
the fact that the lifted action to $S^*Q$ is free, the secondary
and C-L stratifications coincide with the partition
\eqref{semifreequotient}. As for every $(H)\in I_Q$ different from
$(e)$ we have $(H)\in I_Q\setminus I^*_Q$, the contact seams
$CS^{(H)}:=CS_{(H)\succ (e)}$ are Legendrian submanifolds of
$\mathcal{C}_0$. The dimension of each connected component is then
given by formula \eqref{dimensionseam} noting that
$\mathrm{dim\,}Q_{(e)}=\mathrm{dim\,}Q$ and
$\mathrm{dim\,}Q^{(H)}=\mathrm{dim\,}Q_{(H)}-\mathrm{dim\,}G+\mathrm{dim\,}H
=0$ for every $(H)\in I_Q\setminus I^*_Q$, since the action on $Q$
is almost semifree.
 \end{proof}
Recall that a group action is called \textit{semifree} if it is
free everywhere except for a set of isolated fixed points.
Semifree actions are important particular cases of almost semifree
actions and they are commonly found in examples. The following
example explicitly illustrates the geometric constructions of this
paper in that situation.

\noindent \textbf{Example: $S^1$ acting on $S^*\mathbb{R}^2$.}
Consider $Q=\mathbb{R}^2$ with Euclidean coordinates $(x_1,x_2)$
and its cotangent bundle $T^*\mathbb{R}^2=\mathbb{R}^2\times
\mathbb{R}^2$ with coordinates $(x_1,x_2,y_1,y_2)$. The action of
$S^1$ by rotations on $\mathbb{R}^2$ (a semifree action with
${\mathbb{R}^2}_{(S^1)}=\{(0,0)\}$) lifts to $T^*\mathbb{R}^2$ by
the induced diagonal action. A Hilbert basis for the ring of
$S^1$-invariant polynomials for this cotangent lifted action is
given by (see \cite{CuBa1977}, \S1.4)
$$\begin{array}{lll}
\sigma_1 & = & x_1^2+x_2^2+y_1^2+y_2^2,\\
\sigma_2 & = & 2(x_1y_1+x_2y_2),\\
\sigma_3 & = & y_1^2+y_2^2-x_1^2-x_2^2,\\
\sigma_4 & = & x_1y_2-x_2y_1.\end{array}
$$
These polynomials satisfy the semialgebraic relations
$$\sigma_1\geq 0,\quad \sigma_1^2=\sigma_2^2+\sigma_3^2+4\sigma_4^2.$$
We can identify the cosphere bundle $S^*\mathbb{R}^2$ with the
subset of $T^*\mathbb{R}^2$ given by the constraint
$$\sigma_1+\sigma_3=2.$$
The cotangent lifted action restricts to $S^*\mathbb{R}^2$ giving
the free lifted action by contactomorphisms. Its associated
momentum map is given by
$$J(x_1,x_2,y_1,y_2)=\sigma_4$$
for $(x_1,x_2,y_1,y_2)\in S^*\mathbb{R}^2$. Consequently, using
invariant theory, the contact reduced space $J^{-1}(0)/S^1$ is
identified with the semialgebraic variety of
$\mathbb{R}^3=\{\sigma_2,\sigma_3,\sigma_1\}$ defined by
$$\mathcal{C}_0\simeq\left\{(\sigma_2,\sigma_3,\sigma_1)\in\mathbb{R}^3\,
:\, \sigma_1\geq
0,\,\sigma_1^2=\sigma_2^2+\sigma_3^2,\,\sigma_1+\sigma_3=2\right\}.$$
This contact reduced space is in fact a smooth manifold since it
is the parabola obtained intersecting the plane
$P=\{\sigma_1+\sigma_3=2\}$ with the upper half of the cone
$\sigma_1^2=\sigma_2^2+\sigma_3^2$. Its smooth structure is
induced from the ambient space $\mathbb{R}^3$. This was to be
expected since the action on the contact manifold
$S^*\mathbb{R}^2$ is free.

However, this reduced space is no longer a cosphere bundle since
the action on the base is semifree. We investigate now how the
stratified bundle structure of $\mathcal{C}_0$ obtained in the
previous sections arises here. Note that $Q/G=\mathbb{R}^2/S^1$
can be identified with the subset of $\mathbb{R}^3$ given by
$$Q/G=\left\{(0,-t,t)\, :\, t \geq 0\right\},$$ which is a
half-open line parallel to the plane $P$ containing
$\mathcal{C}_0$. According to the notation employed in this
section, $Q/G$ is a stratified space with strata $Q^{(e)}$ and $*
=(0,0,0)$. The continuous fibration $k^0:\mathcal{C}_0\rightarrow
Q/G$ is given by $k^0
(\sigma_2,\sigma_3,\sigma_1)=(0,1-\sigma_1,\sigma_1-1)$. Note that
$(k^0)^{-1}(Q^{(e)})=L\coprod R$ and $(k^0)^{-1}(*)=(0,1,1)$ (see
figure \ref{parabolaBB}), where $\mathcal{C}_0=L\coprod R\coprod
\{(0,1,1)\}$. In addition, recall that $Q^{(e)}\simeq \mathbb{R}$
and that $S^*\mathbb{R}=\mathbb{R}\sqcup \mathbb{R}$.

So $(k^0)^{-1}(Q^{(e)})=L\coprod R$ is diffeomorphic to the cosphere
bundle $S^*Q^{(e)}$. The fiber over a point $(0,-t,t)\in Q^{(e)}$
is the pair of points $(2\sqrt{t},1-t,1+t)$ and
$(-2\sqrt{t},1-t,1+t)$ which lie in $L$ and $R$ respectively.
Finally, the point $(0,1,1)$, the minimum of the parabola
$\mathcal{C}_0$, is the seam $CS_{(S^1)\succ (e)}$ lying in the
boundary of $S^*Q^{(e)}$. Finally, since both $\mathcal{C}_0$ and
$S^*Q^{(e)}$ are one-dimensional, their contact structures are
trivial, due to the fact that the corresponding contact
distributions must be zero-dimensional.

\begin{figure}\label{parabolaBB}
\begin{center}
\includegraphics[scale=0.5]{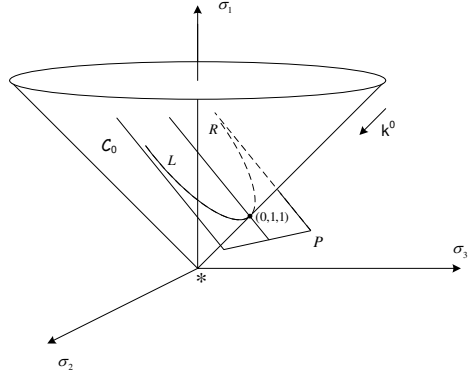}
\end{center}
\caption{The contact reduced space as a parabola fibrating over a
half-closed line.}
\end{figure}

\section{Example: diagonal toral action on $\mathbb{R}^2 \times
\mathbb{R}^2$}

We illustrate the main results obtained in this
paper with one more example rich enough to show all the extra
structure appearing in the cosphere bundle singular reduction.
This time, the reduced contact space $\mathcal{C}_0$ will have
dimension bigger than one and will have hence a non-trivial contact
structure.

Consider the proper action of $G=\TT^2$ on $Q=\RR^2\times\RR^2$,
where each $S^1$ factor acts by rotations on the corresponding
$\RR^2$ factor. The isotropy lattice for this action is shown in
Figure \ref{examples}, where the subconjugation partial order is
represented by arrows. Also, the corresponding stratification
lattice is shown. A stratification lattice is a graphical
arrangement of all the strata of a stratified space where for any
two strata $A,B$ with $A\subseteq\overline{B}$ and such that there
is no other stratum $C$ with the properties $A\subseteq
\overline{C}$ and $C\subseteq \overline{B}$ we write $A\rightarrow
B$. For the action under study, we have $I_Q\setminus I^*_Q=\{
\left(\TT^2\right)\}$.

Let $(x,y)=(x_1,x_2,y_1,y_2)$ be the Euclidean coordinates of a
point in $Q$ and $z=(x,y,u,v)=(x_1,x_2,y_1,y_2,u_1,u_2,v_1,v_2)$
the ones of a covector in $T^*Q\simeq\RR^4\times\RR^4$. The ring of
$G$-invariant polynomials on $T^*Q$ is generated by
$$\begin{array}{lclclcl}
\rho_1 & = & \Vert x\Vert^2+\Vert u\Vert^2&\quad & \sigma_1 & = &
\Vert y\Vert^2+\Vert v\Vert^2 \\
\rho_2 & = & 2(x\cdot u) &\quad & \sigma_2 & = & 2(y\cdot v)  \\
\rho_3 & = & \Vert u\Vert^2-\Vert x\Vert^2 &\quad & \sigma_3 & = &
\Vert v\Vert^2-\Vert y\Vert^2  \\
\rho_4 & = & x_1u_2-x_2u_1 &\quad & \sigma_4 & = & y_1v_2-y_2v_1. \\
& & &
\end{array}
$$
These polynomials, which form a Hilbert basis, are subject to the
following semi-algebraic relations
$$
\rho_1\ge0,\quad\sigma_1\ge0,\quad\rho_1^2=\rho_2^2+\rho_3^2+4\rho_4^2,\quad
\sigma_1^2=\sigma_2^2+\sigma_3^2+4\sigma_4^2.$$

 Identifying the cosphere bundle $S^*\RR^4$ with $\RR^4\times S^3\subset \RR^4\times \RR^4$, where
 $S^3=\{(u,v)\in\RR^2\times\RR^2\, :\, \|u\|^2+\|v\|^2=1\}$,
 it is easy
 to see that
its contact structure is given by the kernel of the restriction of
the Liouville one-form $\theta=udx+vdy$ and that the associated
momentum map $J:S^*\RR^4\rightarrow {\RR^2}$ is given by
$J(x,y,u,v)=(\rho_4,\sigma_4)\in\mathbb{R}^2.$ Consequently, we
still have two more constraints to describe the zero-momentum
level set:
$$\rho_4=0\quad\text{and}\quad\sigma_4=0.$$ Notice that we can
also see $S^*\RR^4$ as the subset of $\RR^8$ defined by the
additional constraint $$\rho_1+\rho_3+\sigma_1+\sigma_3=2.$$

The associated $G$-invariant Hilbert map is defined by
 $$\gamma:J^{-1}(0)\rightarrow
\RR^3\times\RR^3\text{,} \gamma(z)=(\rho_1(z),\rho_2(z),\rho_3(z);
\sigma_1(z), \sigma_2(z),\sigma_3(z)),$$ and we can identify the
reduced contact space with the image of $\gamma$, i.e. with the
semialgebraic variety of $\RR^6$ defined by
$$
\mathcal{C}_0\simeq
\{(\mathbf{\rho};\mathbf{\sigma})\in\RR^6\,:\,\rho_1,\sigma_1\ge0,\,
\rho_1^2=\rho_2^2+\rho_3^2,\,
\sigma_1^2=\sigma_2^2+\sigma_3^2,\,\rho_1+\rho_3+\sigma_1+\sigma_3=2\}
$$
which is the intersection between the product of two cones,
$C_1\times C_2$ and the hypersurface
$H:=\{(\rho_1,\rho_3,\sigma_1,\sigma_3)\in\RR^4\,:\,\rho_1+\rho_3+\sigma_1+\sigma_3=2\}.$
(see Figure \ref{cone-doble}).

 The Reeb vector field on
$S^*\mathbb{R}^4$ is given by $R(x,y,u,v)=(u,v,0,0)$ for any
$(x,y,u,v)\in \RR^4\times S^3$ and the flow of the corresponding
reduced Reeb vector field on $\mathcal{C}_0$ at a point
$(\mathbf{\rho_0};\mathbf{\sigma_0})$ is easily computed as
$$\begin{array}{ll}
\rho_1(t) = \rho_{01}+\rho_{02}t+\frac 12 (\rho_{01}+\rho_{03})t^2\vspace{1mm}\\
\rho_2(t)=\rho_{02}+(\rho_{01}+\rho_{03})t\vspace{1mm}\\
\rho_3(t)=\rho_{03}-\rho_{02}t-\frac 12(\rho_{01}+\rho_{03})t^2\vspace{1mm}\\
\sigma_1(t) = \sigma_{01}+\sigma_{02}t+\frac 12(\sigma_{01}+\sigma_{03})t^2\vspace{1mm}\\
\sigma_2(t)=\sigma_{02}+(\sigma_{01}+\sigma_{03})t\vspace{1mm}\\
\sigma_3(t)=\sigma_{03}-\sigma_{02}t-\frac
12(\sigma_{01}+\sigma_{03})t^2.
\end{array}$$

Applying Proposition \ref{propofirst}, we know that the orbit
types of $J^{-1}(0)$ are exactly those given by $I^*_Q$ and hence
the contact strata of $\mathcal{C}_0$ are in bijective
correspondence with the strata of $Q$ given by $I^*_Q$. We then
have
$$\begin{array}{lll}
T^*Q_{\mathbf{e}}&=&\{(x,y,u,v)\in\RR^8\, :\, (x,y)\neq\mathbf{0}\}\\
T^*Q_{(S^1\times
e)}&=&\{(\mathbf{0},y,\mathbf{0},v)\in\RR^8 \, :\, y\neq\mathbf{0}\}\\
T^*Q_{e\times
S^1}&=&\{(x,\mathbf{0},u,\mathbf{0})\in\RR^8 \, :\, x\neq\mathbf{0}\}\\
N^*Q_{\mathbf{e}}&=&\{(x,y,\mathbf{0},\mathbf{0}) \, :\,
x\neq\mathbf{0}\,
\text{,}\,y\neq\mathbf{0}\}\\
N^*Q_{S^1\times
e}&=&\{(\mathbf{0},y,u,\mathbf{0}) \, :\, y\neq\mathbf{0}\}\\
N^*Q_{e\times S^1}&=&\{(x,\mathbf{0},\mathbf{0},v) \, :\, x\neq\mathbf{0}\}.\\
\end{array}
$$
Consequently, a direct computation gives the following orbit types
for the zero momentum map
$$\begin{array}{lll}
J^{-1}(0) & = & \{ z\in\RR^4\times S^3:\rho_4(z)=\sigma_4(z)=0\}\vspace{1mm} \\

(J^{-1}(0))_{(\mathbf{e})} & = &  \left\{z\in J^{-1}(0): x\neq
0,\, y\neq 0\right\}\\ && \coprod \left\{z\in
J^{-1}(0) \, :\, x=0\,\text{,}\,y\neq  0,\,u\neq 0\right\}\\
& & \coprod \left\{z\in J^{-1}(0) \, :\, x \neq 0, y=0,\,v\neq
0\right\}\\ & & \coprod
\left\{z\in\{0_{\RR^4}\}\times S^3 \, :\, u \neq 0,\,v\neq 0\right\}\vspace{1mm}\\

(J^{-1}(0))_{(e\times S^1)} & = & \left\{z\in J^{-1}(0) \, :\,
y=v=0,\,x\neq 0
 \right\}\\ && \coprod \{z\in J^{-1}(0) \, :\, x=y=v=0,\,\|u\|=1  \}
\vspace{1mm}\\

(J^{-1}(0))_{(S^1\times e)} & = & \left\{z\in J^{-1}(0) \, :\,
x=u=0,\,y\neq 0
 \right\}\\ && \coprod \{z\in J^{-1}(0) \, :\, x=y=u=0,\,\|v\|=1
 \}.
\end{array}
$$
 Using the image of the Hilbert map $\gamma$
we can realize the contact strata given by Theorems \ref{WL1},
\ref{theoremsecondary}, and \ref{coisotropicstratification} as:
$$
\begin{array}{lll}
 \mathcal{C}_0^{(\mathbf{e})} & = &CC_{(\mathbf{e})}\coprod
CS_{(S^1\times e)\succ (\mathbf{e})}\coprod CS_{(e\times S^1)\succ
(\mathbf{e})}\coprod CS_{(\TT^2)\succ
(\mathbf{e})}\vspace{1.5mm}\\

 CC_{(\mathbf{e})}
&=&\left\{(\rho;\sigma): \rho_1,\sigma_1>0
,\,\rho_1\neq\rho_3,\,\sigma_1\neq\sigma_3
,\,\rho_1^2=\rho_2^2+\rho_3^2,\right.\\ & & \left.
\sigma_1^2=\sigma_2^2+\sigma_3^2,\,
\rho_1+\rho_3+\sigma_1+\sigma_3=2\right\}\vspace{1mm}\\

CS_{(S^1\times e)\succ (\mathbf{e})}
&=&\left\{(\rho;\sigma):\rho_1,\sigma_1>0,\,\sigma_1\neq\sigma_3,\,\rho_1=\rho_3,\,
\rho_2=0, \right.\\ & & \left.
2\rho_1+\sigma_1+\sigma_3=2,\,\sigma_1^2=\sigma_2^2+\sigma_3^2\right\}\\
 &=&(\RR_+\times C_2)\cap\left\{2\rho_1+\sigma_1+\sigma_3=2,\,\sigma_1\neq\sigma_3\right\}\vspace{1mm}\\

 CS_{(e\times S^1)\succ
(\mathbf{e})}&=&\left\{(\rho;\sigma):\rho_1,\sigma_1>0,\,\rho_1\neq\rho_3,\,\sigma_1=\sigma_3,\,
\sigma_2=0, \right. \\  & & \left. 2\sigma_1+\rho_1+\rho_3=2,\,\rho_1^2=\rho_2^2+\rho_3^2\right\}\\
 &=&(C_1\times\RR_+)\cap\left\{2\sigma_1+\rho_1+\rho_3=2,\,\rho_1\neq\rho_3\right\}\vspace{1mm}\\

CS_{(\TT^2)\succ
(\mathbf{e})}&=&\left\{(\rho;\sigma):\rho_1,\sigma_1>0,\,\rho_1=\rho_3,\,\sigma_1=\sigma_3,\,
\rho_2=\sigma_2=0, \right. \\ & & \left. \rho_1+\sigma_1=1\right\}\vspace{1mm}\\

 \mathcal{C}_0^{(e\times S^1)} & = &CC_{(e\times S^1)}\coprod
CS_{(\TT^2)\succ (e\times S^1)}\vspace{1mm}\\

CC_{(e\times
S^1)}&=&\left\{(\rho;\mathbf{0}):\rho_1>0,\,\rho_1+\rho_3=2,\,
\rho_1^2=\rho_2^2+\rho_3^2\right\}\setminus \{(1,0,1;\mathbf{0})\}\vspace{1mm}\\

CS_{(\TT^2)\succ (e\times S^1)}&=&\left\{(1,0,1;0,0,0)\right\}\vspace{1mm}\\

 \mathcal{C}_0^{(S^1\times e)} & = & CC_{(S^1\times e)}\coprod
CS_{(\TT^2)\succ (S^1\times e)}\vspace{1mm}\\

CC_{(S^1\times
e)}&=&\left\{(\mathbf{0};\sigma):\sigma_1>0,\,\sigma_1+\sigma_3=2,\,
\sigma_1^2=\sigma_2^2+\sigma_3^2\right\}\setminus \{(\mathbf{0};1,0,1) \}\vspace{1mm}\\

CS_{(\TT^2)\succ (S^1\times e)}&=&\left\{(0,0,0;1,0,1)\right\}.\\
\end{array}
$$

The corresponding contact, secondary and C-L stratification
lattices in $\mathcal{C}_0$ are shown in Figure
\ref{examplesecondary}. Notice that $(\mathbf{e})$ is the
principal orbit type in $Q$ and, therefore, $CC_{(\mathbf{e})}$ is
open and dense in the reduced space $\mathcal{C}_0$. The contact
seams $CS_{(\TT^2)\succ (S^1\times e)}$, $CS_{(\TT^2)\succ (e\times
S^1)}$, and $CS_{(\TT^2)\succ (\mathbf{e})}$ are Legendrian
submanifolds of their contact strata, while the rest are
coisotropic. Every contact seam is mapped by the flow of the
reduced Reeb vector field into the CC-secondary stratum of its
corresponding contact stratum as it can be easily checked.

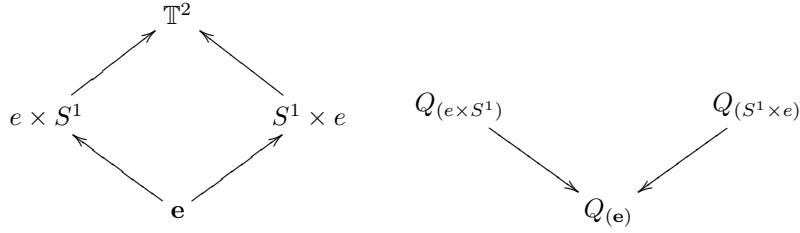
\begin{figure}
\begin{center}
$$\xymatrix{ & \TT^2  &  \\
e\times S^1 \ar[ur] & & S^1\times e \ar[ul]\\
 & \mathbf{e} \ar[ul] \ar[ur] &
}\qquad
\xymatrix{ &   &  \\
Q_{\left(e\times S^1\right)} \ar[dr] & & Q_{(S^1\times e)} \ar[dl]\\
 &Q_{(\mathbf{e})}&
}$$ \caption{\label{examples}Isotropy and stratification lattices
for the $\TT ^2$ action on $\RR ^4$.}
    \end{center}
\end{figure}

\begin{figure}
\begin{center}
\includegraphics[scale=0.7]{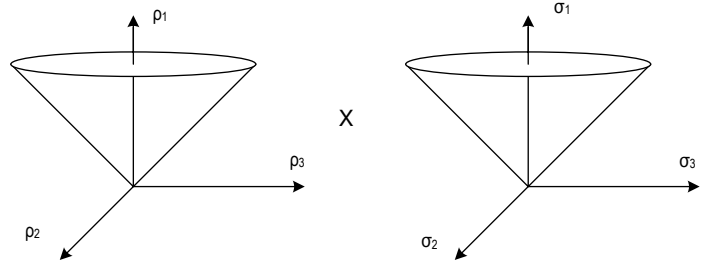}
\end{center}
\caption{The ambient space of $\mathcal{C}_0$} \label{cone-doble}
\end{figure}

\begin{center} 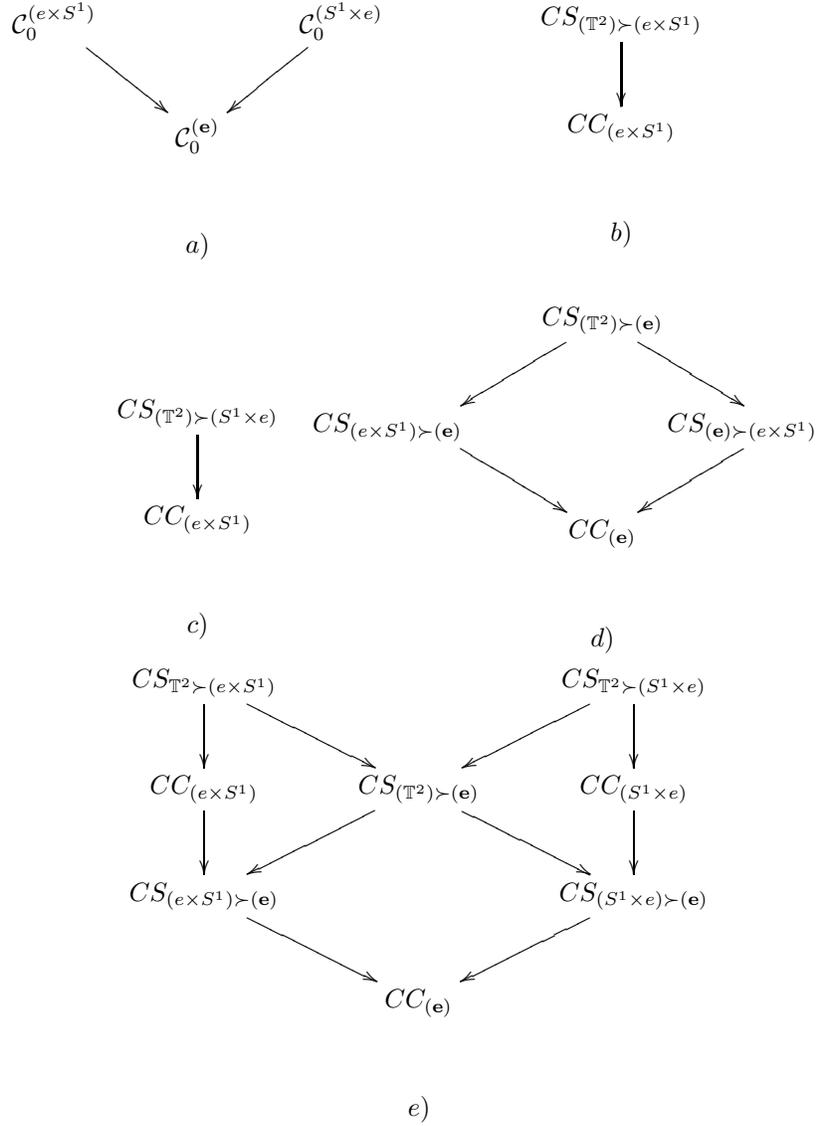
\begin{figure}
   $$\begin{array}{c}   \begin{array}{cc}

\xymatrix{
\mathcal{C}_0^{(e\times S^1)}\ar[dr] & & \mathcal{C}_0^{(S^1\times e)}\ar[dl]\\
  & \mathcal{C}_0^{(\mathbf{e})} &\\ & a) & }

&

\xymatrix{ CS_{(\TT^2)\succ (e\times S^1)}\ar[d]\\ CC_{(e\times
S^1)}\\  b)  }
\\
&
\\

\xymatrix{\\ CS_{(\TT^2)\succ
(S^1\times e)}\ar[d]\\
CC_{(e\times S^1)}\\  c) }

&

\hspace{-0.5cm}\xymatrix{ &  CS_{(\TT^2)\succ (\mathbf{e})}\ar[dl]\ar[dr] &\\
\hspace{-1cm}CS_{(e\times S^1)\succ (\mathbf{e})}\ar[dr] & &
\hspace{-1cm}CS_{(\mathbf{e})\succ (e\times S^1)}\ar[dl]\\
  & CC_{(\mathbf{e})} &\\ & d) & }
\end{array}

\ \\

\xymatrix{  CS_{\TT^2\succ(e\times S^1)}\ar[d]\ar[dr] & &
CS_{\TT^2\succ(S^1 \times e)}\ar[d]\ar[dl]\\
CC_{(e\times S^1)}\ar[d] & CS_{(\TT^2) \succ (\mathbf{e})}\ar[dl]\ar[dr]  & CC_{(S^1\times e)}\ar[d]\\
CS_{(e\times S^1)\succ (\mathbf{e})}\ar[dr] & &
CS_{(S^1\times e)\succ (\mathbf{e})}\ar[dl]\\
& CC_{(\mathbf{e})} & \\ & e) & }
\end{array}$$\caption{\label{examplesecondary} a) Contact
stratification of $\mathcal{C}_0$. Secondary stratifications of:
b) $\mathcal{C}_0^{(e\times S^1)}$, c) $\mathcal{C}_0^{( S^1\times
e)}$ and d) $\mathcal{C}_0^{(\mathbf{e})}$. e)
Coisotropic-Legendrian stratification of $\mathcal{C}_0$.}
              \end{figure}\end{center}

In order to understand the bundle structure of these
stratifications, we embed $Q$ in $T^*Q$ as the zero section and we
identify $Q/G$ with the subset of the image of $\gamma$ given by
$$
Q/G=\left\{(t_1,0,-t_1;t_2,0,-t_2):t_1,t_2\ge
0\right\}\simeq\RR_+\times\RR_+,
$$
a half-plane parallel to $H$. The strata of its orbit
stratification are
$$
\begin{array}{lll}
Q^{(e\times S^1)}&=&\left\{(t_1,0,-t_1;\mathbf{0}):t_1>0\right\}\\
Q^{(S^1\times e)}&=&\left\{(\mathbf{0};t_2,0,-t_2):t_2>0\right\}\\
Q^{(\mathbf{e})}&=&\left\{(t_1,0,-t_1;t_2,;0,-t_2):t_1,t_2>0\right\}
\end{array}
$$
and we obtain that the corresponding cosphere-like strata of
$\mathcal{C}_0$ are diffeomorphic to the cosphere bundles
$$
S^*Q^{(e\times S^1)}\simeq S^*Q^{(S^1\times e)}\simeq\RR\sqcup\RR
\quad\text{and}\quad S^*Q^{(\mathbf{e})}\simeq\RR^2\times S^1.
$$
 The continuous fibration
$k^0:\mathcal{C}_0\rightarrow Q/G$ is given by
$k^0(\rho_1,\rho_2,\rho_3;\sigma_1,\sigma_2,\sigma_3)=(
\rho_1-1,0,1-\rho_1;\sigma_1-1,0,1-\sigma_1).$

 \end{document}